# Using Rolling Circles to Generate Caustic Envelopes Resulting from Reflected Light

**Abstract.** Given any smooth plane curve $\boldsymbol{\alpha}(s)$representing a mirror that reflects light the usual way and any radiant light source at a point in the plane, the reflected light will produce a caustic envelope. For such an envelope, we show that there is an associated curve $\boldsymbol{\beta}(s)$ and a family of circles C(s) that roll on $\boldsymbol{\beta}(s)$ without slipping such that there is a point on each circle that will trace the caustic envelope as the circles roll. For a given curve $\boldsymbol{\alpha}(s)$ and for all radiants at infinity there is a single curve $\boldsymbol{\beta}(s)$ and family of circles C(s) that roll on $\boldsymbol{\beta}(s)$ so that the different points on C(s) will simultaneously trace out, as the circles roll, all caustic envelopes from these radiants at infinity. We explore many classical examples using this method.

## 1. Introduction

A commonly observed phenomenon is the pattern on the bottom of a cylindrical container resulting from reflected light off of its vertical sides (see photo 1). This brightly lit piecewise smooth curve with a single cusp, sometimes called the "coffee cup caustic" is the envelope of reflected light rays. Each point on the curve is a focal point for some point on the reflective surface, with the focal point depending on the location of the light source and the curvature of the surface.

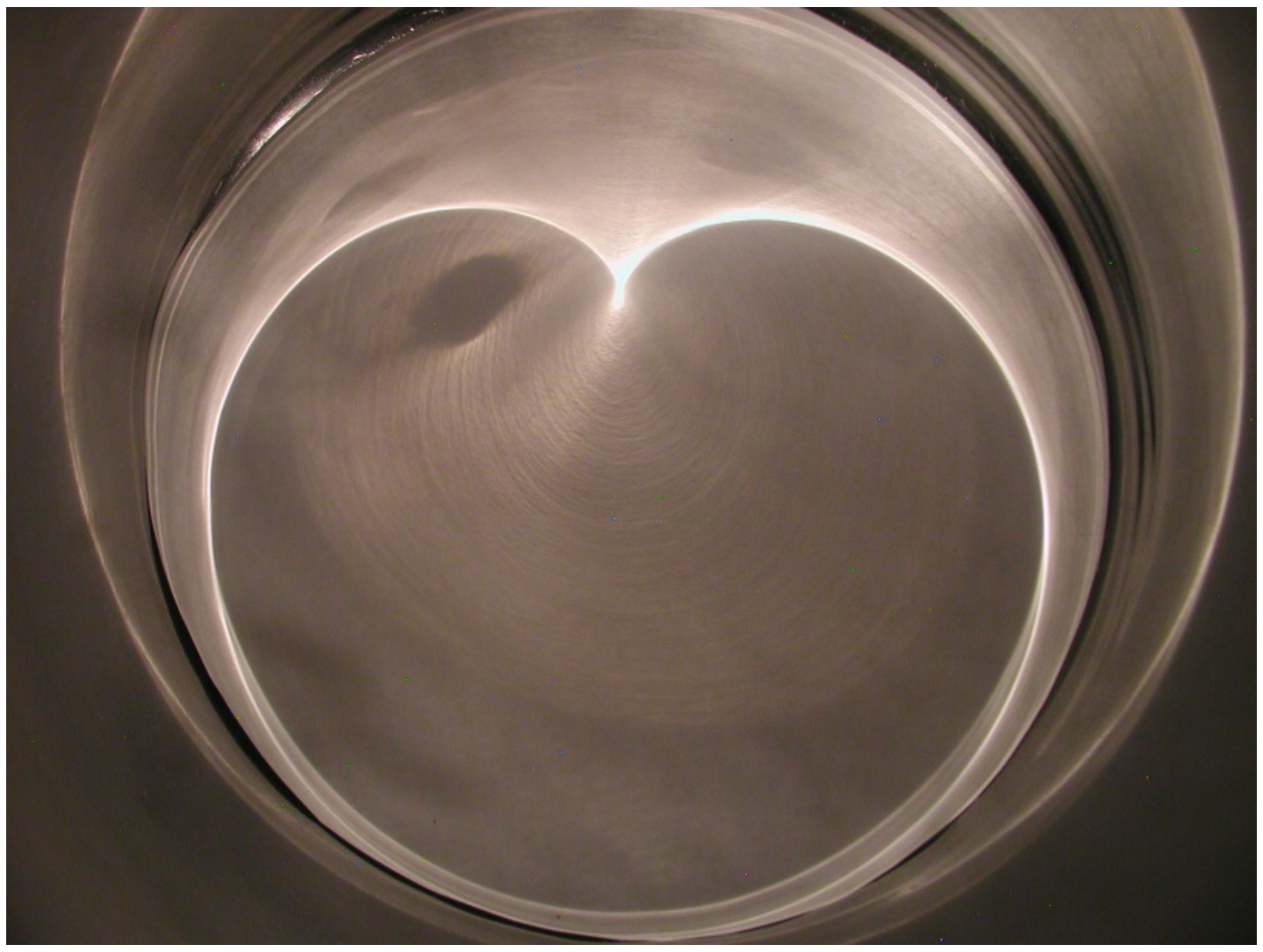

Photo 1

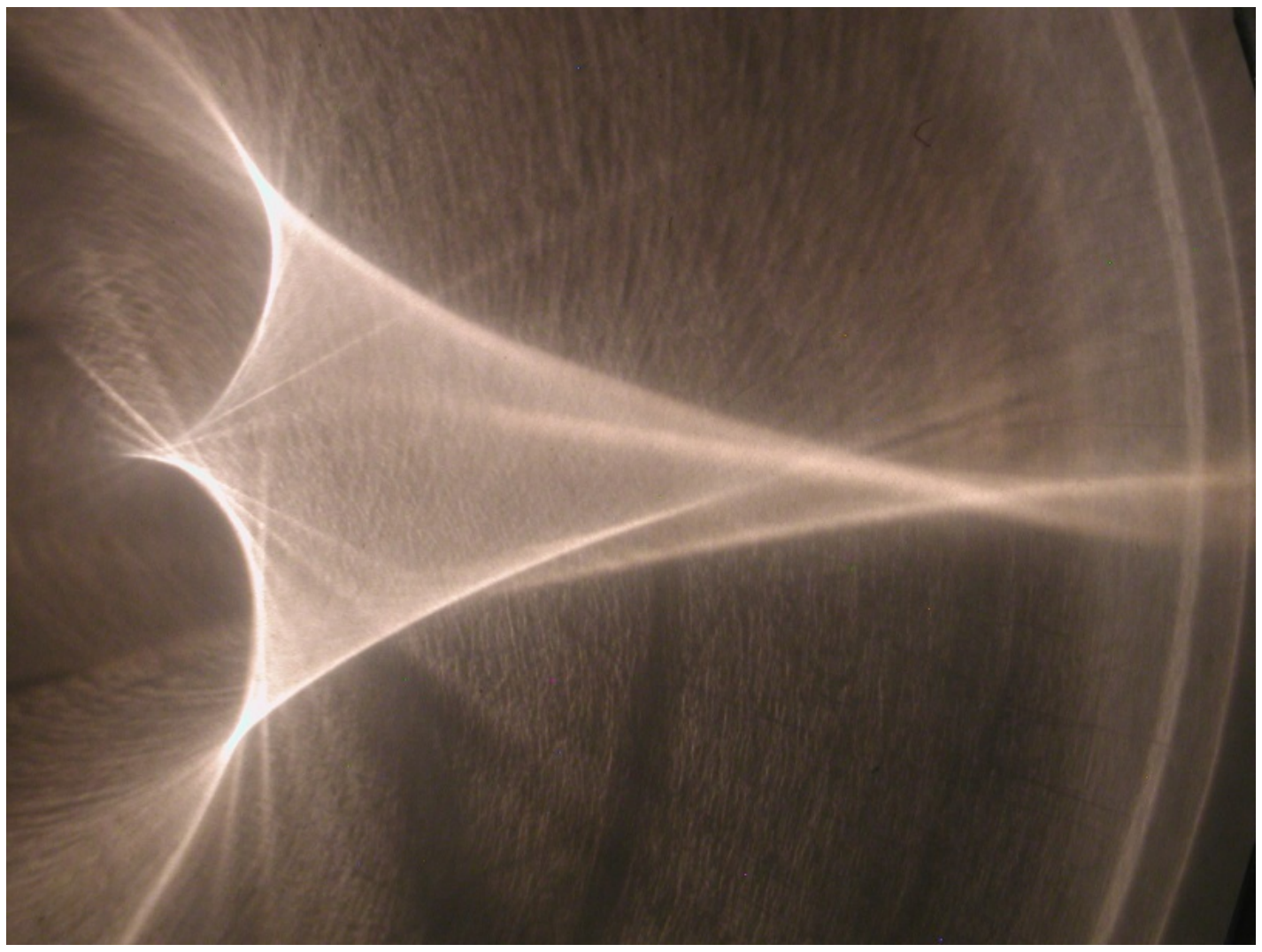

Photo 2

Different reflective surfaces can produce a variety of different caustic curves and even a fixed reflective surface may produce an infinite number of caustics as the position of the light source is changed. Photos 1 and 2 are of caustics created by shining a flash light into a deep stock pot. In the first photo the source is near the rim, and in the second the source is in the interior but off center. Figure 1 shows two caustics resulting from internal reflection in an elliptical mirror with the light sources at the dots (each caustic being concentrated in the opposite quadrant from the source). The reflected rays are shown but the incident rays are omitted.

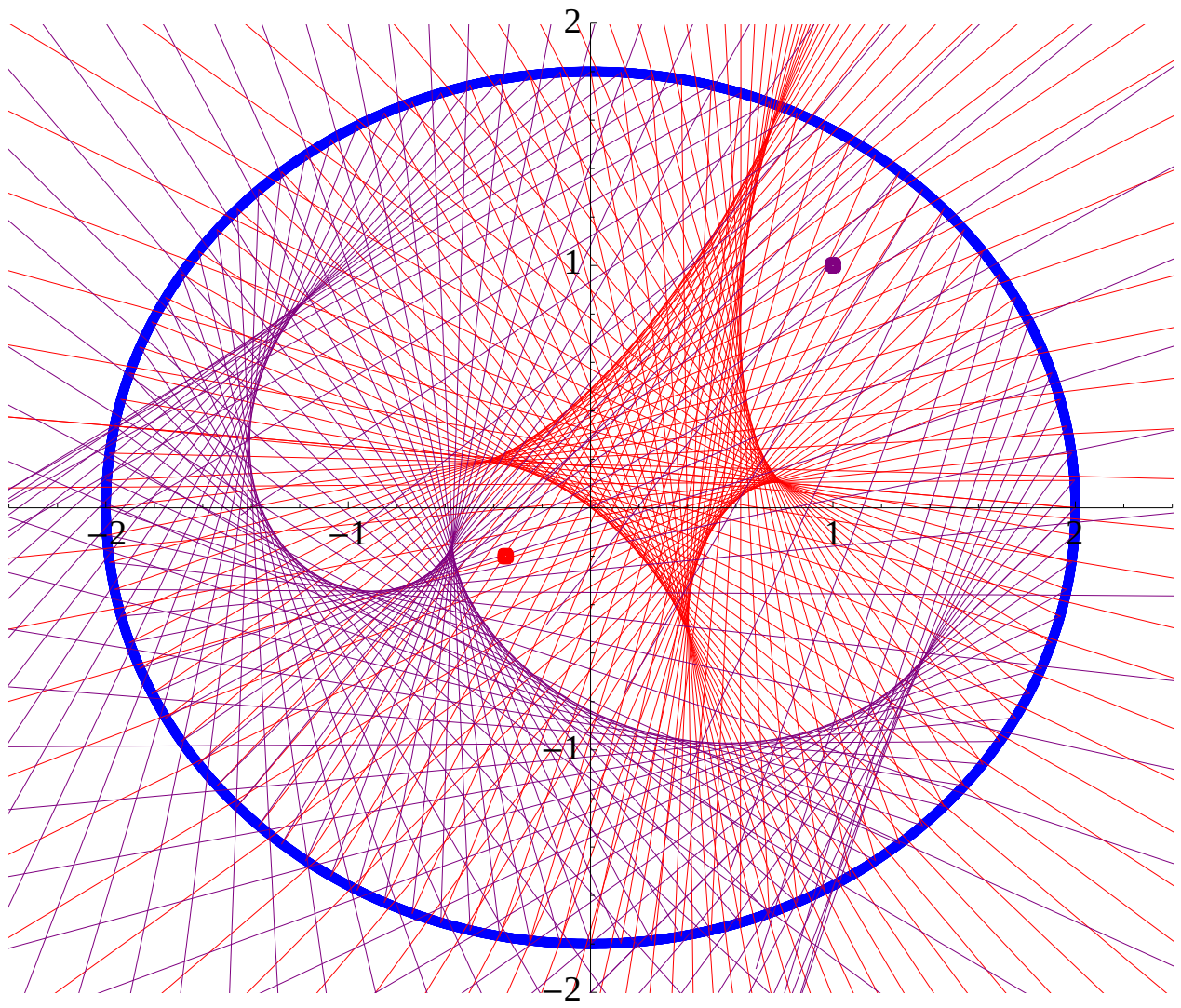

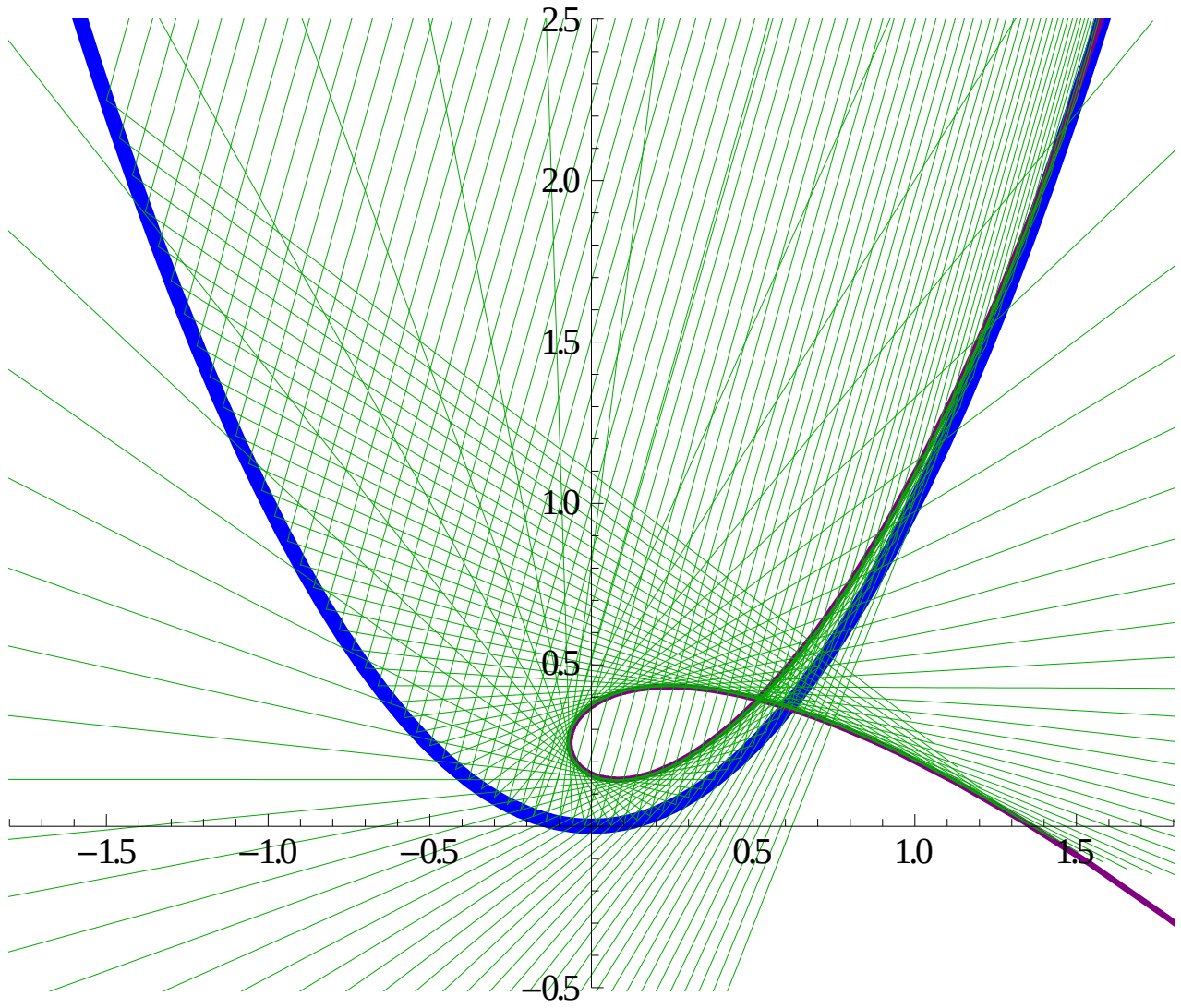

**Figure 1** Two caustics from internal reflection in an elliptical mirror

**Figure 2** Caustic from a radiant at infinity in a parabolic mirror

Figure 2 shows a caustic from a parabolic mirror with light source at infinity in a direction 75° counter clockwise from the positive x-axis.

If the object that produces the caustic is a surface, then the caustic itself is generally a surface too. In the case of the coffee cup caustic, the bottom of the container is a screen that produces a planar cross section of the caustic. In this article, we will confine ourselves to 2 dimensions, where light from a point source (called a *radiant*) in a plane (or at infinity) is reflected off a smooth curve in that plane. There is a long history of the study of caustics. A. Cayley [11] wrote a memoir in 1857 describing caustics both from reflections and refractions, giving calculations for all cases involving a circular mirror (see [4] for equations). See also [8, 12, 16, and 19]. Lockwood [18] "A Book of Curves" discusses caustics, gives a table of common curves and some of their associated caustics and provides a useful glossary for plane curves. The optical caustics in this article share much with the theory of mathematical billiards. See for example [15 or 23].

## The Coffee Cup Caustic

In the article [19], the authors derive a parameterization for the caustic in the case where the reflective curve is the standard unit circle $\boldsymbol{\alpha}(\theta) = \langle \cos\theta, \sin\theta \rangle$. Here the light source is at infinity in the direction of the negative x-axis with light traveling in a direction parallel to the positive x-axis. The caustic envelope **E** has parameterization

$$\mathbf{E}(\theta) = \left\langle \tfrac{3}{4}\cos\theta - \tfrac{1}{4}\cos 3\theta,\ \tfrac{3}{4}\sin\theta - \tfrac{1}{4}\sin 3\theta \right\rangle = \tfrac{3}{4}\langle \cos\theta, \sin\theta \rangle - \tfrac{1}{4}\langle \cos 3\theta, \sin 3\theta \rangle.$$

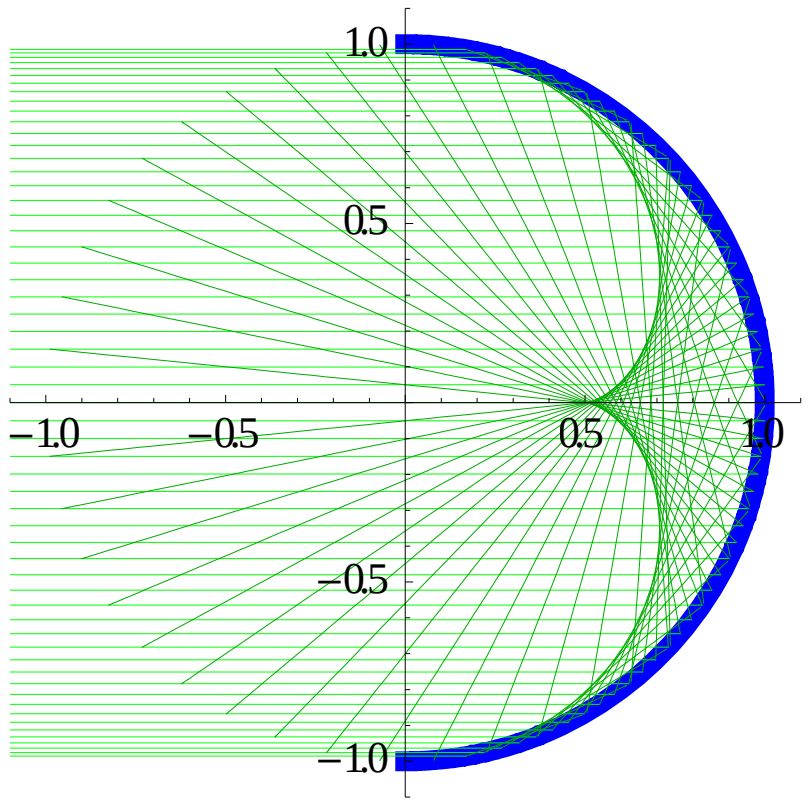

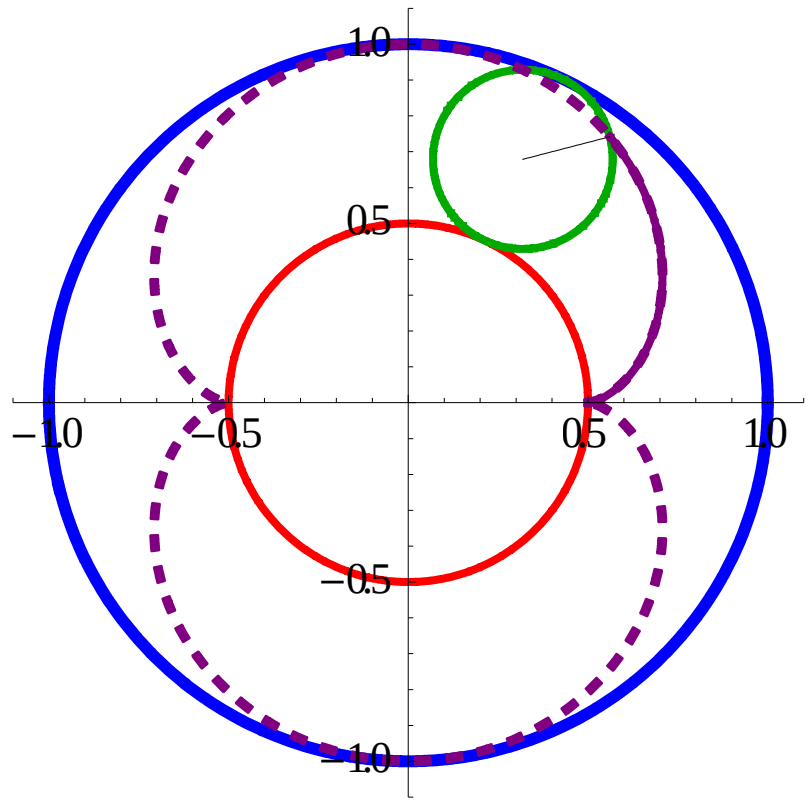

**Figure 3** Light reflecting in a semi-circular mirror

**Figure 4** The caustic as an epicycloid

The second version of **E**($\theta$) as a sum of two vector functions reveals a curious alternative description of the caustic as an epicycloid traced out by a point on a circle of radius ¼ rolling around the outside of a circle of radius ½.

The main theorem of this article is that every caustic can be generated in this manner.

**Theorem 1** Given a caustic envelope **E** resulting from a given reflective curve $\boldsymbol{\alpha}$(s) and radiant S (for source), there is a curve $\boldsymbol{\beta}(s)$ and a family of circles $C_s$ that roll on $\boldsymbol{\beta}$ without slipping such that each circle has a point that will trace the caustic envelope as the circles roll.

To see animations demonstrating Theorem 1 for a variety of different curves $\boldsymbol{\alpha}$, consult [1] for radiants at infinity and [2] for finite radiants. These demonstrations allow the user to select the direction of the radiant at infinity. The radiant direction may also be animated. For more images of caustics see [3].

The caustic is generated somewhat in the manner of a roulette, although not a true roulette because the rolling circles do not generally have a fixed diameter.

In Figure 5, the $\boldsymbol{\alpha}$ curve is the outer ellipse and $\boldsymbol{\beta}$ is the inner elliptical-like curve. The radiant is at infinity with light coming in parallel to the arrow above. The circle rolls around $\boldsymbol{\beta}$ while remaining tangent to $\boldsymbol{\alpha}$. The point on the circle traces out the caustic as the circle rolls.

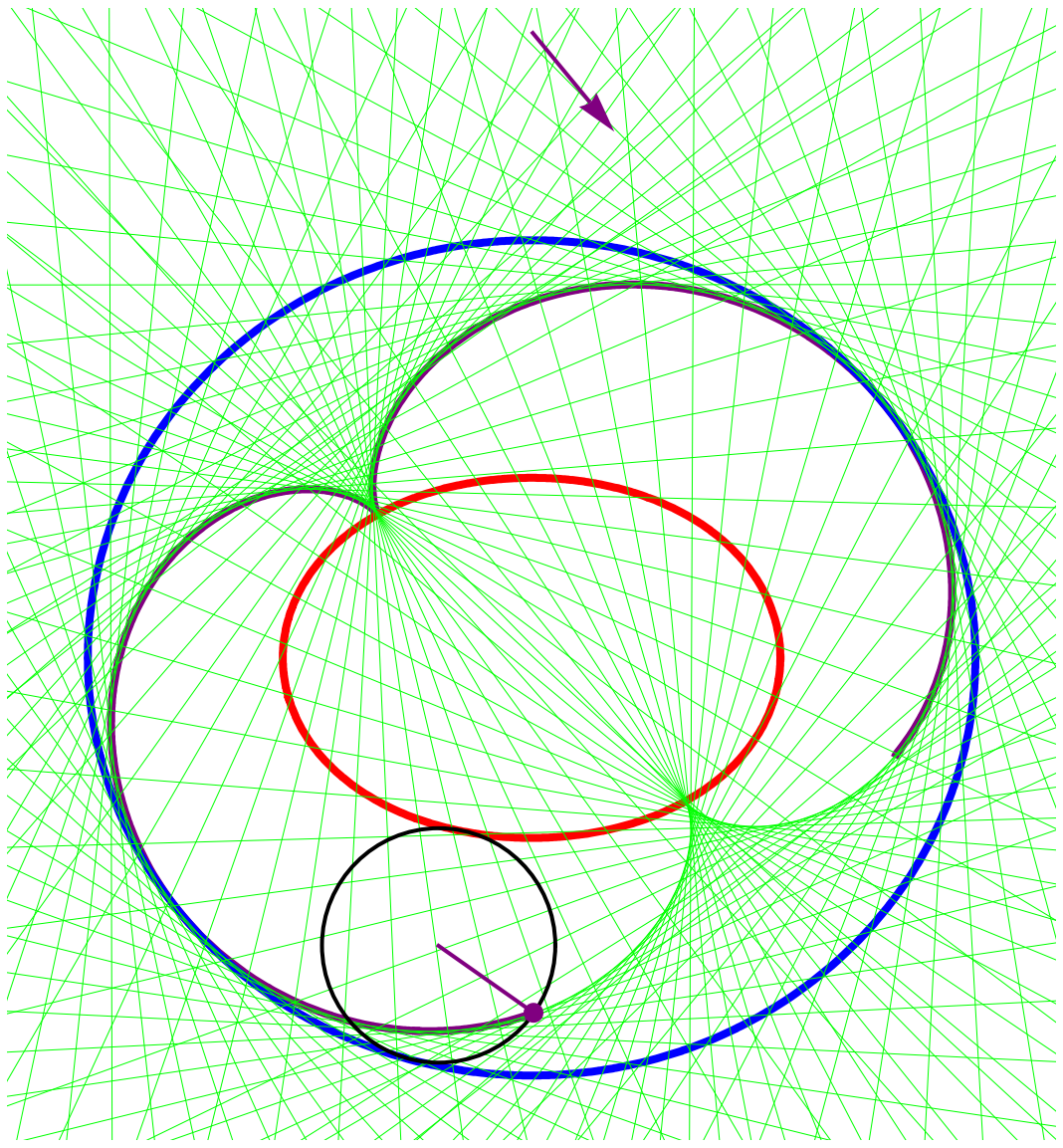

**Figure 5** Illustrating Theorem 1 for an elliptical mirror and radiant at infinity

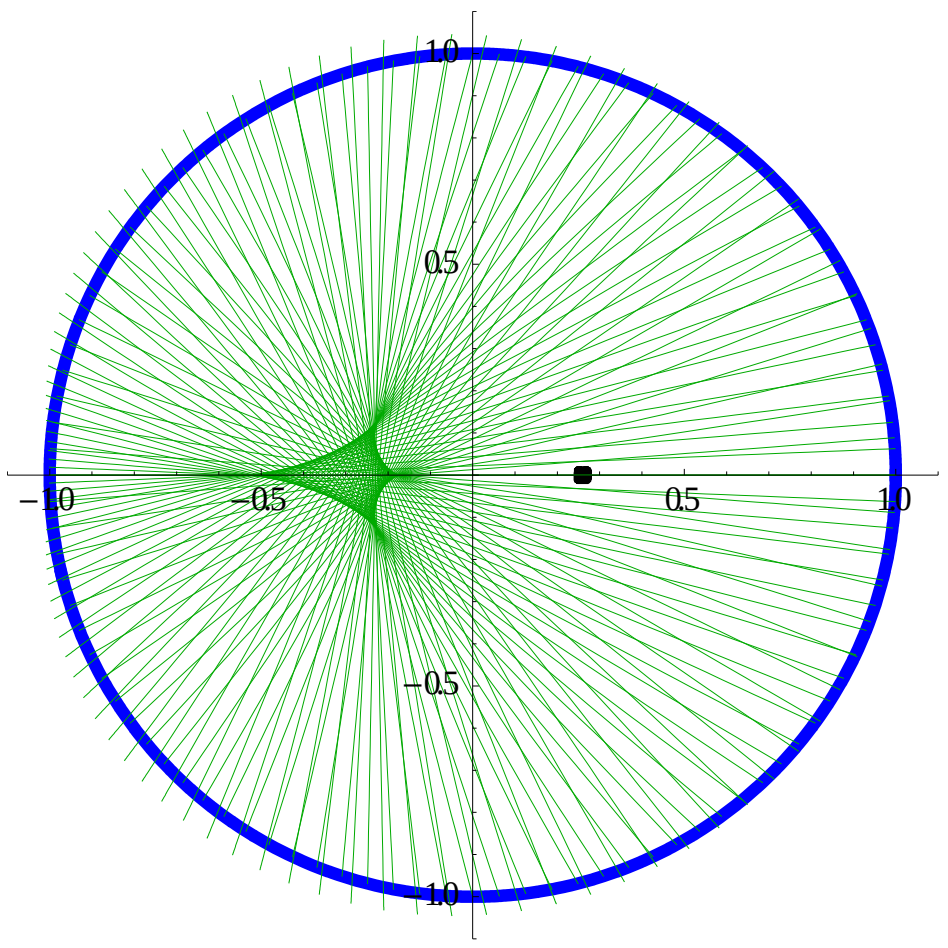

**Figure 6** Internal reflection circular mirror

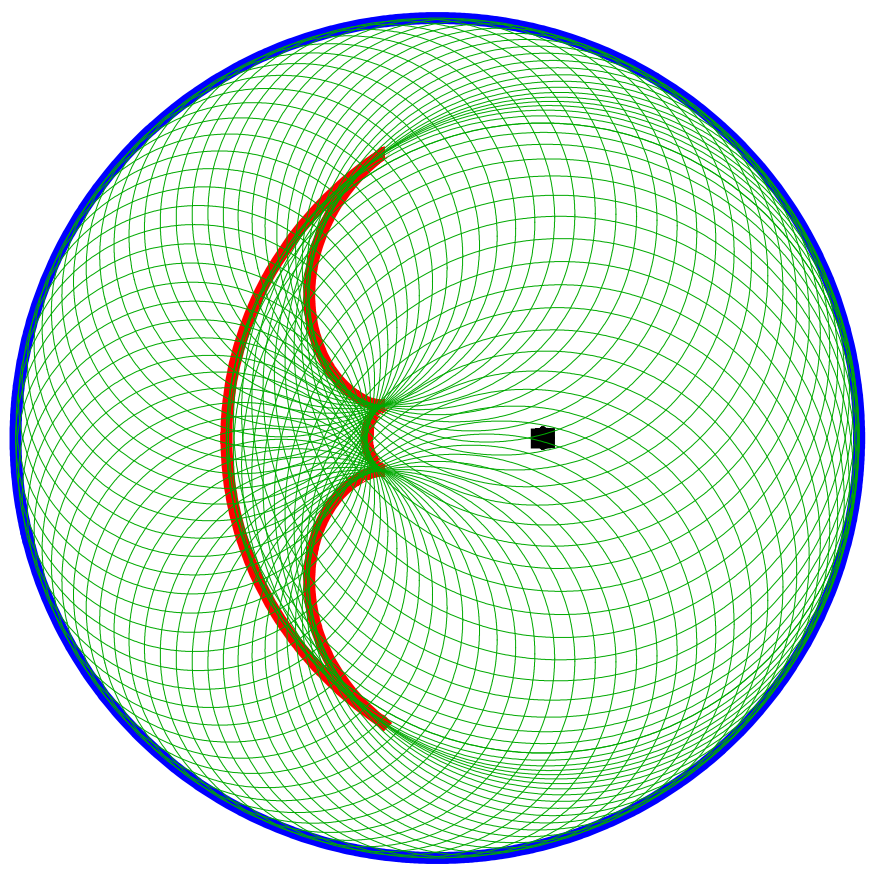

**Figure 7** Circles $C_s$ and $\boldsymbol{\beta}$

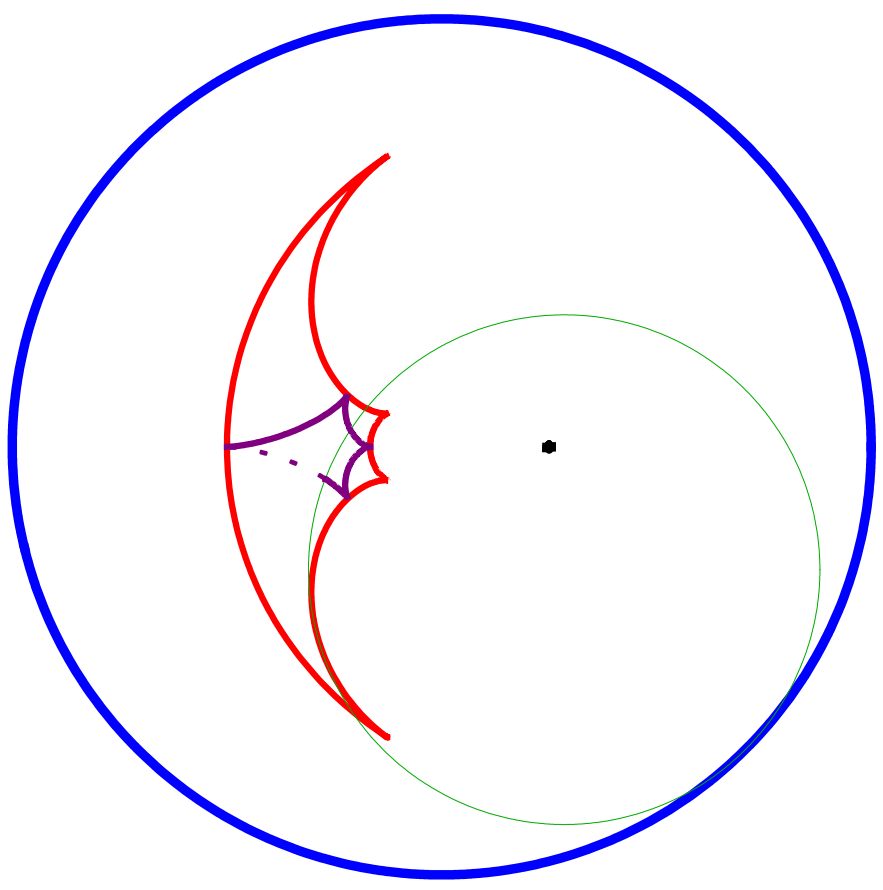

**Figure 8** Tracing the caustic

Figure 6 illustrates internal reflection in a circular mirror from a radiant at the dot creating the apparent caustic envelope. Figure 7 shows the associated family of circles and their envelope $\boldsymbol{\beta}$, and Figure 8 is a "stop action" shot of the circles rolling on $\boldsymbol{\beta}$ with a point tracing out the caustic envelope.

Throughout this paper, $\boldsymbol{\alpha}$(s) will represent a smooth curve that will reflect light rays in the usual way – the angle of incidence equals the angle of reflection. The parameter s is arc length along the curve $\boldsymbol{\alpha}$. However, on any other curve associated with $\boldsymbol{\alpha}$ the same parameter s will not be arc length. The curve $\boldsymbol{\alpha}$ has curvature $\kappa = \kappa(s)$ and radius of curvature r = r(s) = 1/ $\kappa(s)$.

## 2. Focal Circles

Let S be a radiant for $\boldsymbol{\alpha}$. Light from S reflected from $\boldsymbol{\alpha}$ at $\boldsymbol{\alpha}$(s) and $\boldsymbol{\alpha}(s+\delta)$ will intersect at a point **F**($\delta$). We say light from radiant S reflected at $\boldsymbol{\alpha}$(s) will focus at

$$\mathbf{F} = \lim_{\delta\to 0} \boldsymbol{F}(\delta),$$

provided the limit exists.

The point F is the point of tangency of the reflected ray with the caustic envelope. Equation 1 in Theorem 2 below is the classical equation from geometric optics (see, for example, [20 or 23]). However, we state the mirror equation in a form that is absolutely essential for everything in this paper, and provide a short general proof that requires no special assumptions about how a curve reflects light other than the usual angle of incidence equals the angle of reflection.

**Theorem 2** **(Generalized Mirror Equation)** Let $C_1$ and $C_2$ be circles tangent to $\boldsymbol{\alpha}$ at $\boldsymbol{\alpha}$(s) and on the inside of $\boldsymbol{\alpha}$ at $\boldsymbol{\alpha}$(s) (their centers are in the direction of the normal to the curve) with diameters $d_1$ and $d_2$. Light from a radiant on one circle reflected at $\boldsymbol{\alpha}$(s) will have a focus on the other circle if

$$\frac{1}{d_1}+\frac{1}{d_2}=\frac{2}{r}=2\kappa \qquad (1)$$

where $\kappa$ is the curvature and $r=\frac{1}{\kappa}$ is the radius of curvature. In other words, the harmonic mean of $d_1$ and $d_2$ is the radius of curvature.

To prove Theorems 1 and 2 we will repeatedly utilize the following basic lemma which tells us the rotation rate of a line segment based on the motion of the endpoints.

**Lemma 1** Let $\boldsymbol{u}(t)$ and $\boldsymbol{w}(t)$ be smooth curves and let $\boldsymbol{c}(t) = \boldsymbol{w}(t) - \boldsymbol{u}(t)$ be the vector from $\boldsymbol{u}\ to\ \boldsymbol{w}$. Let $\sigma(t)$ be the direction angle for $\boldsymbol{c}(t)$. Then

$$\frac{d\sigma}{dt} = \frac{\|\boldsymbol{u}'\| cos\ \varphi_1 - \|\boldsymbol{w}'\| cos\ \varphi_2}{d}$$

where $\varphi_1$ and $\varphi_2$ are the angles between $\boldsymbol{c}(t)$ and the normals to $\boldsymbol{u}(t)$ and $\boldsymbol{w}(t)$, respectively, and d is the length of $\boldsymbol{c}(t)$. In particular, if $\boldsymbol{u} = \boldsymbol{u}(s)$ is parameterized by arc length s, and either

(1) $\boldsymbol{w}(t) = S$, a fixed point (the radiant); or

(2) $\boldsymbol{c}(t)$ is tangent to $\boldsymbol{w}(t)$,

then

$$\frac{d\sigma}{ds} = \frac{cos\ \varphi_1}{d}$$

**Proof** Using subscripts to indicate component functions we have $\tan\sigma(t) = c_y(t)/c_x(t)$ and thus

$$\sigma'(t) = \frac{c_x c_y' - c_y c_x'}{\left(1 + c_y^2/c_x^2\right)(c_x^2)} = \frac{c_x\left(w_y' - u_y'\right) - c_y(w_x' - u_x')}{c_x^2 + c_y^2}$$

$$= \frac{< c_x, c_y > \cdot < -u_y', u_x' > - < c_x, c_y > \cdot < -w_y', w_x' >}{d^2} = \frac{d\|\boldsymbol{u}'\| cos\ \varphi_1 - d\|\boldsymbol{w}'\| cos\ \varphi_2}{d^2}$$

For the two special cases we have $\|\boldsymbol{u}'\| = 1$, and either (1) $\|\boldsymbol{w}'\| = 0$, or (2) $\varphi_2 = \frac{\pi}{2}$. Both cases result in the simplified equation 2. ■

**Proof of Theorem 2** See Figure 9. Let $\sigma_1, \sigma_2$ and $\sigma_3$ be direction angles for $\overrightarrow{PS}$, $\vec{N}$ and $\overrightarrow{PF}$, respectively. An elementary argument shows $\frac{d\sigma_1}{ds} = \frac{\cos\varphi}{D_1}$ and $\frac{d\sigma_3}{ds} = \frac{\cos\varphi}{D_2}$. Also we have $\frac{d\sigma_2}{ds} = \kappa$. Since $\sigma_1 + \sigma_3 = 2\sigma_2$,

$$\frac{d\sigma_1}{ds} + \frac{d\sigma_3}{ds} = 2\frac{d\sigma_2}{ds} \quad \text{or} \quad \frac{\cos\varphi}{D_1} + \frac{\cos\varphi}{D_2} = 2\kappa \qquad (2)$$

Define $d_1$ and $d_2$ by

$$D_1 = d_1 \cos\varphi \ \text{ and } \ D_2 = d_2 \cos\varphi \qquad (3)$$

Substituting into (2) gives the equation $\frac{1}{d_1} + \frac{1}{d_2} = 2\kappa$.

Points $(D_1, \varphi)$ and $(D_2, \varphi)$ satisfying (3) describe, in polar coordinates, the circles defined in Theorem 2. ■

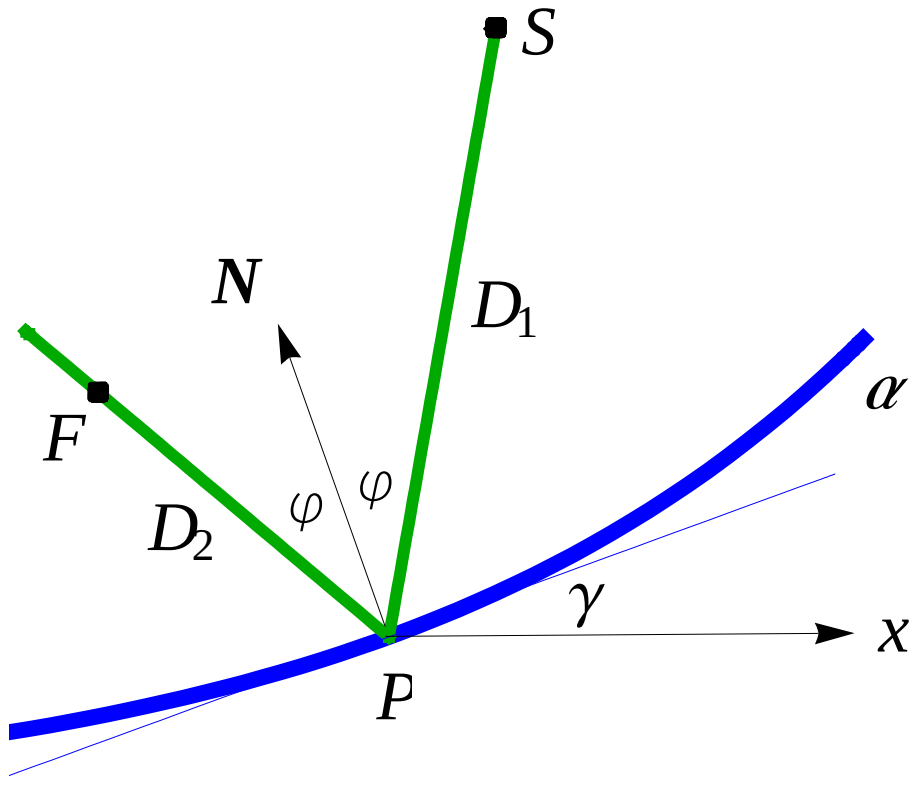

**Figure 9** Angles and distances for proof of Theorem 2

The animations at [1] can be used to illustrate Theorem 2 in the case where the radiant is at infinity. Selecting the point P of reflection and animating the direction of the radiant, the point of the caustic associated with P will move around the associated focal circle for P.

Let $C_d$ be the circle of diameter d tangent to $\boldsymbol{\alpha}$ at P and let r be the radius of curvature at P. The circle $C_{r/2}$ is called the *discriminant circle* (see [8]). There are a number of simple consequences of Theorem 2 that we now list.

(i) If the radiant is outside $C_r$ , the focus is on the inside of $C_r$, and vice versa.
(ii) If the radiant is on $C_r$ , the focus will be on $C_r$ .
(iii) If the radiant is at infinity, then the focus is on the discriminant circle $C_{r/2}$ (and vice versa).
(iv) If the radiant were inside the discriminant circle, then reflected light will diverge rather than focus.
(v) No light can focus inside the discriminant circle.

If the reflected light rays diverge, the standard procedure is to extend the reflected ray to an entire line running both directions from the point of reflection. The proof of Theorem 2 works equally well in this case with the following modification. If the incident angle is $\varphi$, then the angle of "reflection" properly measured is $\varphi + \pi$, and then either $d_1$ or $d_2$ would be negative. The corresponding circle tangent at the point of reflection is on the outside (convex side) of the curve. We can reinterpret properties (iv) and (v) listed above as follows.

(iv) If the radiant is inside the discriminant circle, then the reflected line has focus on the outside of the curve.
(v) If the radiant is outside the curve, then the focus is inside the discriminant circle.

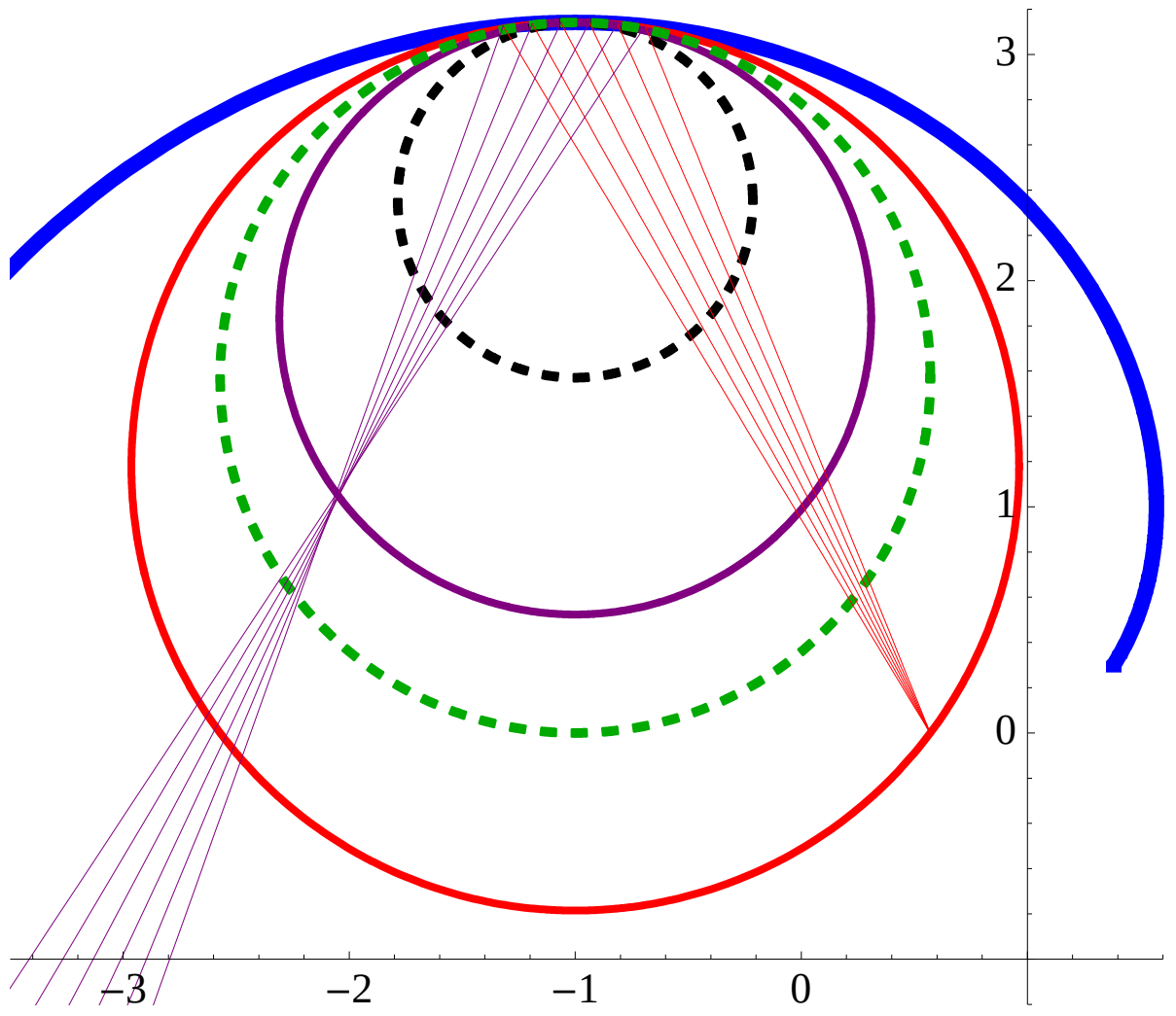

**Figure 10** Any radiant on the outer solid circle will focus on the inner solid circle.

Figure 10 illustrates the theorem with light reflected from the spiral curve near a point where the radius of curvature is $\pi$. The smaller dotted circle is the discriminant circle with diameter one half of the radius of curvature and the larger dotted circle has diameter equal to the radius of curvature. Light from any radiant on the outer solid circle of diameter $\frac{5\pi}{4}$ will focus at a point on the smaller solid circle of diameter $\frac{5\pi}{6}$.

In [8], geometrical methods are used to study the location of the focus F as a result of reflection at P = $\boldsymbol{\alpha}$ $(s_0)$ from the source at S. It is shown that there is a unique conic that has at least three point contact with $\boldsymbol{\alpha}$ (s) at P such that S is at one of the foci of the conic. The point F will then be located at the other focus of the conic. Many interesting properties of the caustic envelope are derived from this method. The unique conic is an ellipse, parabola, or hyperbola depending on whether the source S is outside, on, or inside the discriminant circle, respectively.

## 3. The Envelope Theorem

Let $\boldsymbol{\alpha}$ be a smooth curve and S be a radiant (possibly at infinity). Then for each point $\boldsymbol{\alpha}$(s) of $\boldsymbol{\alpha}$, there is a circle $C_s$ tangent to $\boldsymbol{\alpha}$ at $\boldsymbol{\alpha}(s)$ such that the light from S reflected at $\boldsymbol{\alpha}(s)$ will focus at some point of $C_s$ as described by Theorem 2. Call $C_s$ the *focal circle* for $\boldsymbol{\alpha}$ at $\boldsymbol{\alpha}(s)$ (relative to S). The curve $\boldsymbol{\alpha}$ is then an envelope of the family of circles $C_s$. There is a second envelope for the family $C_s$, which we denote by $\boldsymbol{\beta}$(s) (see Figure 11). As we shall see, generally the curve $\boldsymbol{\beta}$ is complicated to parameterize. It should be observed that the line segment joining $\boldsymbol{\alpha}(s)$ and $\boldsymbol{\beta}$(s) is generally not a diameter for the focal circle $C_s$. The location of $\boldsymbol{\beta}$(s) on $C_s$ relative to $\boldsymbol{\alpha}(s)$ depends on the relative rate of change in the curvature at $\boldsymbol{\alpha}(s)$ and the location of the radiant S.

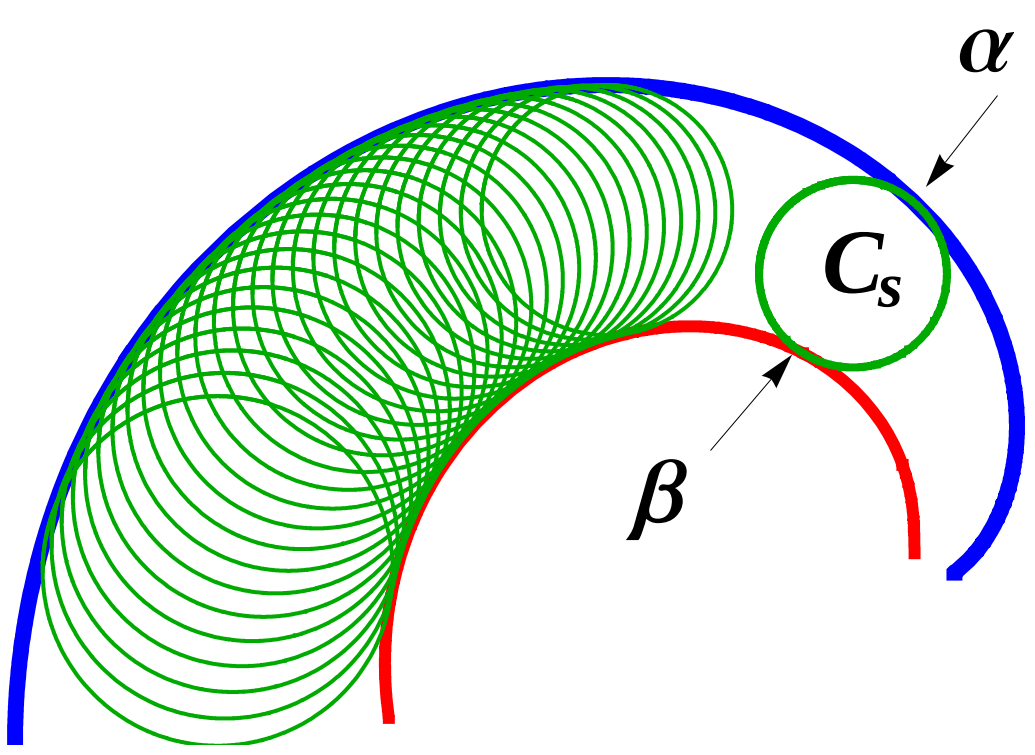

**Figure 11** Focal circles and the two envelopes

**Theorem 3 (The Envelope Theorem)** Let $\boldsymbol{\alpha}(s)$ = <$\alpha_x, \alpha_y$> be a smooth curve parameterized by arc length s, with nonzero curvature $\kappa(s)$. Let $C_s$ be a circle tangent to $\boldsymbol{\alpha}$ at $\boldsymbol{\alpha}(s)$, with $C_s$ having center in the direction of the normal to $\boldsymbol{\alpha}(s)$. Furthermore, let $C_s$ have radius R(s), which is a differentiable function of s. Then the family of circles determines a second envelope $\boldsymbol{\beta}$(s) given by

$$\boldsymbol{\beta}(s) = \boldsymbol{\alpha}(s) + R(s)\big(1+\cos\big(2\delta(s)\big)\big)\boldsymbol{N}(s) + R(s)\sin\big(2\delta(s)\big)\,\boldsymbol{T}(s) \qquad (4)$$

where $\tan(\delta(s)) = \frac{R'(s)}{R(s)\kappa(s)-1}$ , $\boldsymbol{T}(s) =< \alpha_x'(s), \alpha_y'(s) >$, and $\boldsymbol{N}(s) = \boldsymbol{T}'(s)/\|\boldsymbol{T}'(s)\|$

are the unit tangent and normal vectors, respectively.

**Remarks:** For simplicity, Theorem 2 is stated in the case $\kappa(s) \neq 0$. However, Theorem 2 can be extended so this restriction can be omitted by selecting a continuous unit normal vector and also allowing the circles to be on the convex side. Also for convenience, the angle $\delta(s)$ has been defined via $\tan(\delta(s))$. If $R(s_0)\kappa(s_0) - 1 = 0$, then the tangent is not defined but the angle $\delta(s_0)$ is defined ($\pm\frac{\pi}{2}$), and Equation 4 still applies. In this case, $C_{s_0}$ is the osculating circle at $\boldsymbol{\alpha}(s_0)$, and $\boldsymbol{\beta}(s_0) = \boldsymbol{\alpha}(s_0)$. Hence, $\boldsymbol{\alpha}(s)$ and $\boldsymbol{\beta}(s)$ intersect whenever $C_s$ is an osculating circle (assuming $R'(s) \neq 0$).

**Proof** The circle $C_s$ has center $\boldsymbol{\alpha}(s) + R(s)\,\boldsymbol{N}(s)$ and radius R(s). Therefore the family of circles $C_s$ can be described by the family of equations

$$(\boldsymbol{x} - \alpha_x(\boldsymbol{s}) + R(\boldsymbol{s})\alpha_y'(s))^2 + (\boldsymbol{y} - \alpha_y(\boldsymbol{s}) - R(\boldsymbol{s})\alpha_x'(s))^2 - R^2(\boldsymbol{s}) = 0$$

which has the form F(x,y,s) = 0. The standard method to find the envelope(s) of the family is to solve the system

$$\frac{\partial F}{\partial s}(x,y,s) = 0, \quad F(x,y,s) = 0.$$

We leave the details as an exercise for the interested reader, or see [7] for detailed proofs of statements in this paper.

----------

Here it is.

Write $F = P^2 + Q^2 - R^2$ for convenience and compute

$$\frac{\partial F}{\partial s} = 2P\left(-\alpha_x'(s) + R'(s)\alpha_y'(s) + R(s)\alpha_y''(s)\right) + 2Q\left(-\alpha_y'(s) - R'(s)\alpha_x'(s) - R(s)\alpha_x''(s)\right) - 2R(s)R'(s)$$

To locate $\boldsymbol{\beta}(s_0)$ relative to $\boldsymbol{\alpha}(s_0)$ we will use local coordinates centered at $\boldsymbol{\alpha}(s_0)$ with orthonormal basis {$\mathbf{T}(s_0), \boldsymbol{N}(s_0)$}. Hence

$$\boldsymbol{\alpha}(s_0) =< 0,0 >,\ \boldsymbol{\alpha}'(s_0) =< 1,0 >,\ \boldsymbol{\alpha}''(s_0) =< 0, \kappa(s_0) >,\ P(x,y,s_0) = x,\ \ Q(x,y,s_0) = y - R(s_0)$$

Setting $\frac{\partial F}{\partial s}(x,y,s_0) = 0$ (and suppressing $s_0$) we get

$$x(-1+R\kappa) + (y-R)(-R') - RR' = 0 \qquad \text{or} \qquad (R\kappa - 1)x = R'y$$

One solution is <x,y> = <0,0> = $\boldsymbol{\alpha}(s_0)$, confirming that $\boldsymbol{\alpha}(s)$ is one of the envelopes. For the other envelope define the angle $\delta$ by $\tan\delta = \frac{R'}{R\kappa-1}$. Let c be the chord connecting $\boldsymbol{\alpha}(s_0)$ and $\boldsymbol{\beta}(s_0)$. The angle from $\boldsymbol{N}(s_0)$ to c is $\delta$, which we will refer to as the *chord angle*. With a little trigonometry we have

$$c = 2R\cos\delta \qquad x = R\sin 2\delta \qquad y = R(1 + \cos 2\delta)$$

and then $\boldsymbol{\beta}(s_0) =< R\sin 2\delta, R(1 + \cos 2\delta) >$. Dropping the local coordinates we then have

$$\boldsymbol{\beta}(s) = \boldsymbol{\alpha}(s) + R(s)\big(1 + \cos(2\delta(s))\big)\boldsymbol{N}(s) + R(s)\sin\big(2\delta(s)\big)\boldsymbol{T}(s) \qquad \blacksquare$$

---------

Let c be the chord that connects $\boldsymbol{\alpha}(s)$ and $\boldsymbol{\beta}(s)$ in the focal circle $C_s$ and let d be the diameter of $C_s$ that has $\boldsymbol{\alpha}(s)$ as one of its endpoints. The angle between c and d is the angle $\delta(s)$ that appears in Equation 4 which we call the *chord angle*.

The Envelope Theorem applies to any family of circles $C_s$ tangent to the smooth curve $\boldsymbol{\alpha}(s)$, not just to the family of focal circles, assuming the radius function R(s) is differentiable. Here are two such applications, the first is unsurprising, but the second less so.

**Corollary 1** If the circles $C_s$ all have the same radius, then $\delta = 0$ and $\boldsymbol{\alpha}(s)$ and $\boldsymbol{\beta}(s)$ are endpoints of a diameter of $C_s$ orthogonal to $\boldsymbol{\alpha}(s)$.

Thus the curves $\boldsymbol{\alpha}(s)$ and $\boldsymbol{\beta}(s)$ would then be so called Bertrand curves, sharing the same normal line for each s.

**Corollary 2** If the circles $C_s$ are the family of osculating circles for $\boldsymbol{\alpha}(s)$, then there is only one envelope ($\boldsymbol{\alpha}(s)$).

**Proof** For osculating circles, R(s) = $\frac{1}{\kappa(s)}$ making $\tan\delta$ undefined. This would put $\boldsymbol{\beta}(s)$ at a right angle to the normal $\boldsymbol{N}(s)$, but on $C_s$. Hence $\boldsymbol{\beta}(s) = \boldsymbol{\alpha}(s)$.

Corollary 2 is essentially the same as an old result known as the Tait-Kneser Theorem ([17, 24]). One could argue that in fact for the family of osculating circles there in no envelope at all, since by Tait-Kneser the osculating circles along an arc with monotone positive curvature are pairwise disjoint and nested. Ordinarily, a point of an envelope of curves should be the limit of intersection points of one family member with a second as the second approaches the first. See also [14] for more about Tait-Kneser and its generalizations (and see the surprising pictures illustrating theTait-Kneser Theorem).

## 4. Radiant at Infinity

In Theorem 1, the focal circles' radii R(s) are determined by the curvature of $\boldsymbol{\alpha}(s)$ and the location of the radiant S as described by Theorem 2. We will often consider the case where the radiant S is at infinity, in which case, by Theorem 2,

$R(s) = \frac{1}{4\kappa(s)}$ , $R'(s) = -\frac{\kappa'(s)}{4\kappa^2(s)}$, and the focal circles are the family of discriminant circles for $\boldsymbol{\alpha}(s)$ . We compute $\tan\delta = \frac{R'}{R\kappa - 1} = \frac{\kappa'(s)}{3\kappa^2(s)}$.

For any smooth curve $\boldsymbol{\alpha}(s)$ with $\kappa(s) \neq 0$ , the quantity a(s) = $-\frac{\kappa'(s)}{3\kappa^2(s)}$ , is called the *aberrancy* of the curve $\boldsymbol{\alpha}(s)$. See [21]. For the current application, we see a(s) = $-\tan\delta = \tan(-\delta)$. Manipulating trigonometric functions, we have

$$\cos 2\delta = \frac{1-a^2}{1+a^2} \qquad \text{and} \qquad \sin 2\delta = \frac{-2a}{1+a^2}\,. \tag{5}$$

Combining equations (4) and (5) yields the following corollary.

**Corollary 3** Let $\boldsymbol{\alpha}(s)$ be a smooth curve with $\kappa(s) \neq 0$, and let $C_s$ be the family of focal circles for $\boldsymbol{\alpha}(s)$ resulting from any radiant at infinity. Then the second envelope $\boldsymbol{\beta}(s)$ is given by

$$\boldsymbol{\beta}(s) = \boldsymbol{\alpha}(s) + \frac{1}{2\kappa(s)(1+a(s)^2)}(\boldsymbol{N}(s) - a(s)\,\boldsymbol{T}(s)), \tag{6}$$

where **N**(s) and **T**(s) are the unit normal and tangent vectors at $\boldsymbol{\alpha}(s)$, respectively.

The concept of the aberrancy of a curve was defined in 1841 by Abel Transon [25] in an attempt to create a geometric or visual quantity that would represent the value of the third derivative of a function y = f(x). For such a curve, Transon defined the aberrancy a and the *angle of aberrancy* $\delta$ by the equation

$$a = \tan\delta = y' - \frac{(1+y'^2)}{3y''^2}y'''$$

Geometrically, at a point P of a curve, let R and Q be points on the curve very near to P such that the chord RQ is parallel to the tangent at P. Let M be the midpoint of such a chord and then form the line through P and M. As the chord approaches the tangent line as a limit, generally the line through P and M will approach a limiting line which is called the *axis of aberrancy*. The angle that the axis of aberrancy makes with the normal at P is the angle of aberrancy $\delta$ described above. Hence, the aberrancy at P is a measure of the curves local departure from symmetry about its normal direction. For much more about aberrancy see the articles by Schot [21, 22] where he also discusses the history of aberrancy. Schot mentions that aberrancy predates even Transon, going back as far as the early 1800's [10]. See also [6, 9, and 26]. The effectiveness of aberrancy as a visual representative of the third derivative is open to question. However, one can verify that the aberrancy $a(s) = -\frac{\kappa'(s)}{3\kappa^2(s)}$ , where the derivative is with respect to arc length s, and hence $|a(s)|$ is intrinsic to the curve (independent of parameterization).

For focal circles with the radiant at infinity, the chord angle $\delta$ is the negative of the angle of aberrancy $\delta$. Thus, the reflection of the axis of aberrancy across the normal at $\boldsymbol{\alpha}(s)$ will contain the point $\boldsymbol{\beta}(s)$.

## 5. Some Examples

**Example 1:** (The Coffee Cup Caustic) The curve $\boldsymbol{\alpha}(s) =< \cos s\,, \sin s >$ has constant curvature $\kappa = 1$, so $R(s) = \frac{1}{4}$ and a = 0. It follows that $\boldsymbol{\beta}(s) =< \cos s\,, \sin s > + \frac{1}{2} < -\cos s, -\sin s > = \frac{1}{2} < \cos s\,, \sin s >$, agreeing with the earlier description of the caustic as an epicycloid traced out by a point on a circle of radius $\frac{1}{4}$ rolling on the outside of a circle of radius $\frac{1}{2}$.

**Example 2:** The spiral $\boldsymbol{\alpha}(t) = < \cos t + t\sin t, \sin t - t\cos t >$ has aberrancy a(t) = $\frac{1}{3t}$ and

$$\boldsymbol{\beta}(t) = \frac{1}{2(9t^2+1)} < (15t^2+2)\cos t + (9t^3+2t)\sin t, (15t^2+2)\sin t + (-9t^3-2t)\cos t >$$

(see Figure 11).

**Example 3:** The parabola $\boldsymbol{\alpha}(t) = < t, bt^2 >$ has aberrancy a(t) = 2bt and $\boldsymbol{\beta}(t) = < 0, \frac{1}{4b} >$, the focus of the parabola.

## 6. Proof of Main Theorem (Outline)

For full details of the Proof, see [7]. We define several direction angles, measured from the positive x-direction, all of which are functions of s:

$\sigma_1$ is the direction angle from $P = \boldsymbol{\alpha}(s)$ to the radiant S,
$\omega$ is the direction from the center of the focal circle $C_s$ to the point to the point of the caustic **E**(s) on $C_s$,
$\gamma$ is the direction angle of the tangent vector $\boldsymbol{\alpha}'(s)$,
$\varphi$ is the angle of incidence/reflection of the light ray,
$d_1$ is the diameter of the circle C tangent to $\boldsymbol{\alpha}(s)$ at P that contains the radiant S.

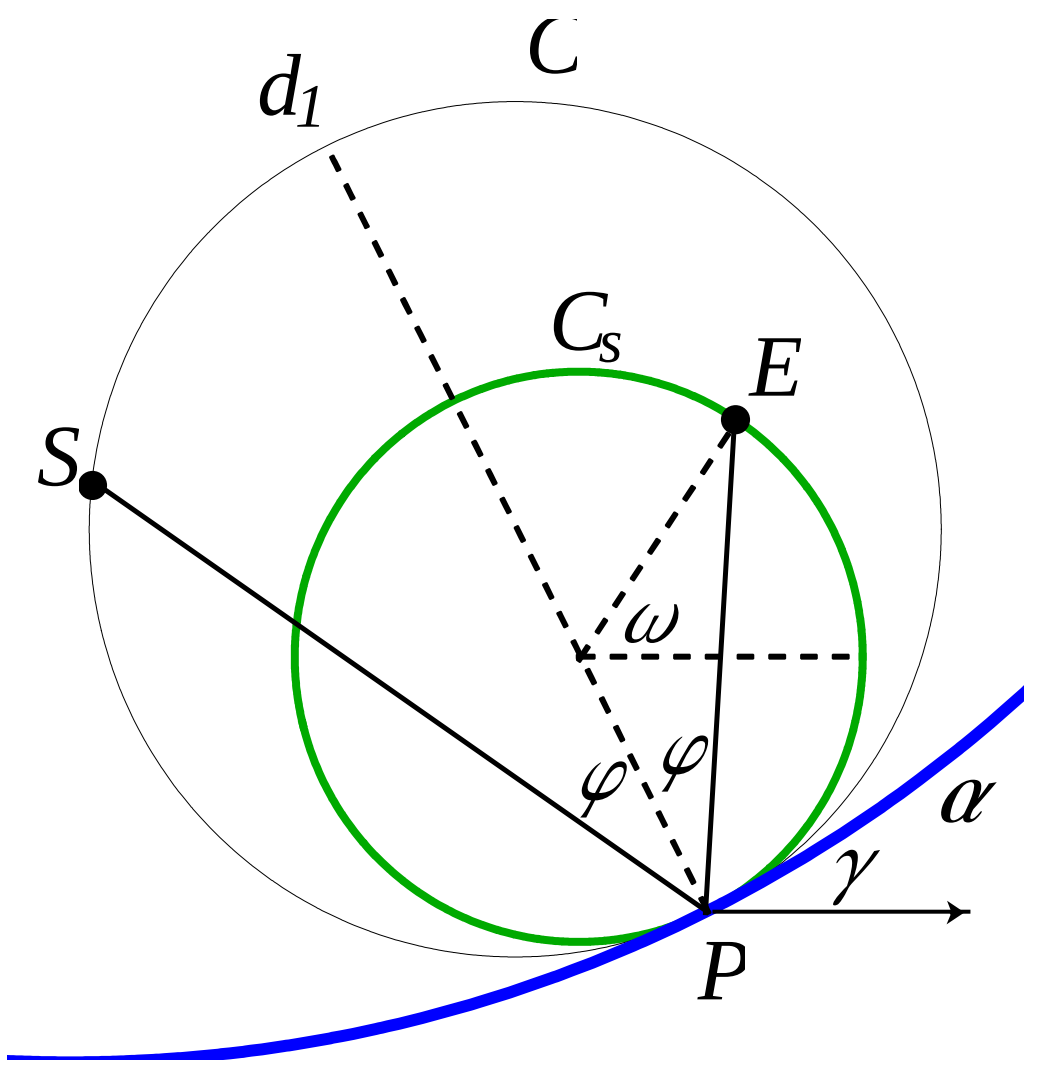

**Figure 12** Definition of the angles

We have $\omega = \frac{3\pi}{2} + 3\gamma - 2\sigma_1$, and as mentioned earlier, $\frac{d\sigma_1}{ds} = \frac{\cos\varphi}{d} = \frac{1}{d_1}$.

From this we compute that the angular velocity of **E**(s) around $C_s$ is $\frac{d\omega}{ds} = 3\frac{d\gamma}{ds} - 2\frac{d\sigma_1}{ds} = 3\kappa(s) - \frac{2}{d_1}$. This angular velocity is independent of the location of **E**(s) around $C_s$. This is an important observation. Every point on $C_s$ is a point of some caustic resulting from reflection at P from a radiant located on the circle C. This "circle of caustics" is, at least momentarily, rotating with uniform angular speed around $C_s$ in the manner of a rigidly rotating wheel.

The point $\boldsymbol{\beta}(s)$ on $C_s$ is also a point of a caustic resulting from reflection at P from some radiant S on C. Fix $s_0$ and let **E**(s) be the caustic from a radiant S such that **E**($s_0$) = $\boldsymbol{\beta}(s_0)$. The proof of Theorem 2 will be complete if we can show the velocity $\mathbf{E}'(s_0) = \mathbf{0}$, for then the motion of **E**(s) around $C_s$ is not only rigid rotation, but also non-slipping in its contact with $\boldsymbol{\beta}(s)$. Note that **E**(s) = $\boldsymbol{\alpha}(s) + R(s) < -\alpha_y'(s), \alpha_x'(s) > + R(s) < \cos\omega, \sin\omega >$ and that the caustic envelope **E**(s) touches $\boldsymbol{\beta}(s)$ if $\omega - \gamma = \frac{\pi}{2} - 2\delta$, where $\delta$ is the chord angle in $C_s$ . Using these observations and using Equation 1, the verification that $\mathbf{E}'(s_0) = \mathbf{0}$ is a straight forward, but long calculation.

-------

Here is the calculation. Write **E**(s) = <x(s), y(s)>, where

x(s) = $\alpha_x(s) \; - R(s)\alpha_y'(s) + R(s)\cos\omega$, y(s) = $\alpha_y(s) \; + R(s)\alpha_x'(s) + R(s)\sin\omega$

We compute

$$x' = \alpha_x' - R'\alpha_y' - R\alpha_y'' + R'\cos\omega - R\sin\omega \cdot \frac{d\omega}{ds}$$

$$y' = \alpha_y' + R'\alpha_x' + R\alpha_x'' + R'\sin\omega + R\cos\omega \cdot \frac{d\omega}{ds}$$

Now fix a point $\boldsymbol{\alpha}(s_0)$ on $\boldsymbol{\alpha}$. Using local coordinates again we have $\boldsymbol{\alpha}(s_0) = <0,0>$, $\boldsymbol{\alpha}'(s_0) = <1,0>$, $\boldsymbol{\alpha}''(s_0) = <0, \kappa(s_0)>$, and so $\gamma = 0$. At $s_0$ the above equations are reduced to

$$x' = 1 - R\kappa + R' \cos\omega - R\sin\omega \cdot \frac{d\omega}{ds}$$

$$y' = R' + R' \sin\omega + R\cos\omega \cdot \frac{d\omega}{ds}$$

The caustic envelopes $\boldsymbol{E}$(s) touches $\boldsymbol{\beta}(s)$ if $\omega - \gamma = \frac{\pi}{2} - 2\delta$, where $\delta$ is the chord angle in $C_s$ and $\gamma$ is the direction angle for the tangent direction $\mathbf{T}(s_0)$. See Figure 12.5.

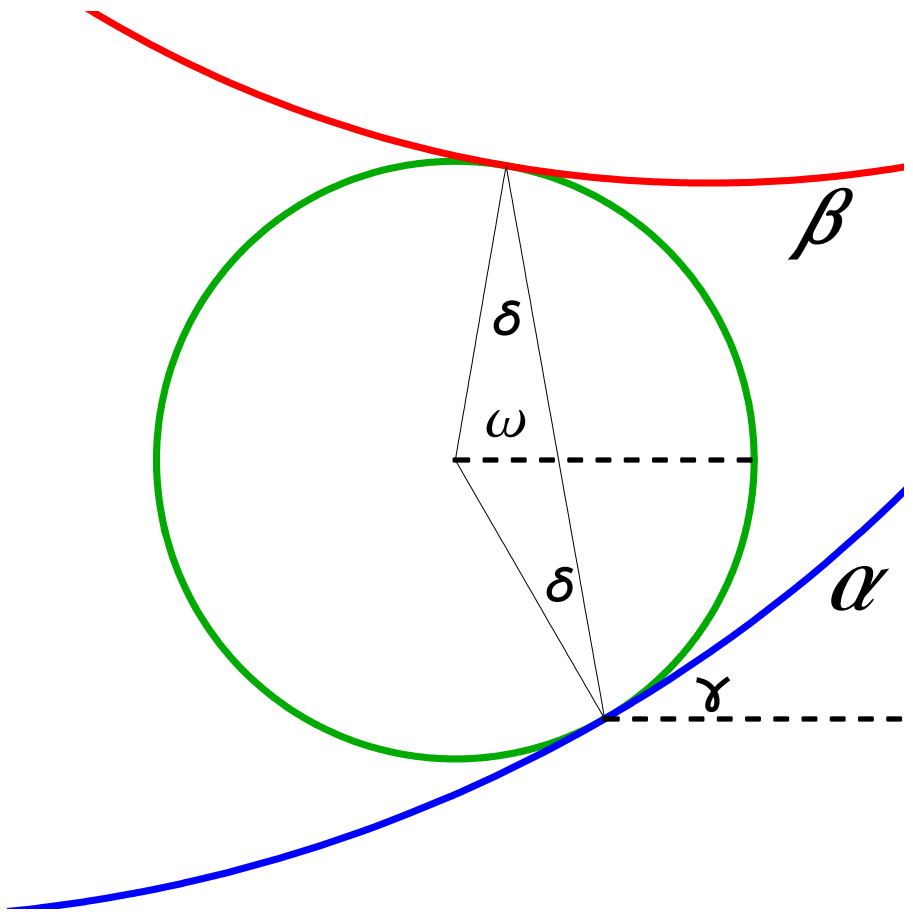

**Figure 12.5** The caustic touches $\boldsymbol{\beta}$

At $\boldsymbol{\alpha}(s_0)$ since $\gamma = 0$, $\omega = \frac{\pi}{2} - 2\delta$. From $\tan\delta = \frac{R'}{R\kappa - 1}$ and some trigonometry we get

$$\cos\omega = \cos(\frac{\pi}{2} - 2\delta) = \sin 2\delta = \frac{2R'(R\kappa - 1)}{(R\kappa - 1)^2 + R'^2}, \qquad \sin\omega = \sin(\frac{\pi}{2} - 2\delta) = \cos 2\delta = \frac{(R\kappa - 1)^2 - R'^2}{(R\kappa - 1)^2 + R'^2}$$

Substituting these into the equation for $x'$ above we get

$$x' = (1 - R\kappa) + R'\left(\frac{2R'(R\kappa - 1)}{(R\kappa - 1)^2 + R'^2}\right) - R\sin\omega \cdot \frac{d\omega}{ds}$$

$$= \frac{(1 - R\kappa)[(R\kappa - 1)^2 + R'^2] + 2R'^2(R\kappa - 1)}{(R\kappa - 1)^2 + R'^2} - R\sin\omega \cdot \frac{d\omega}{ds}$$

$$= \frac{(1 - R\kappa)[(R\kappa - 1)^2 + R'^2 - 2R'^2]}{(R\kappa - 1)^2 + R'^2} - R\sin\omega \cdot \frac{d\omega}{ds}$$

$$= \frac{(1 - R\kappa)[(R\kappa - 1)^2 - R'^2]}{(R\kappa - 1)^2 + R'^2} - R\sin\omega \cdot \frac{d\omega}{ds}$$

$$= (1 - R\kappa)\sin\omega - R\sin\omega \cdot \frac{d\omega}{ds} = \sin\omega\left[1 - R\kappa - R\frac{d\omega}{ds}\right]$$

The radius of the focal circle $C_{s_0}$ is $R = \frac{d_2}{2}$ and by Theorem 2, $\frac{1}{d_1} + \frac{1}{d_2} = 2\kappa$. Combining these we have $R = \frac{d_1}{4\kappa d_1 - 2}$. Also recall that $\frac{d\omega}{ds} = 3\kappa(s_0) - \frac{2}{d_1}$. Substituting these into the above equation for $x'$ yields

$$x' = \sin\omega\,[1 - R\kappa - R(3\kappa - \frac{2}{d_1})] = \sin\omega\ [1 - 4R\kappa + \frac{2R}{d_1}]$$

$$= \sin\omega\left[1 - \frac{4\kappa d_1}{4\kappa d_1 - 2} + \frac{2}{4\kappa d_1 - 2}\right] = \sin\omega \cdot 0 = 0$$

Now we will do the same for $y'$.

$$y' = R'(1 + \sin\omega) + R\cos\omega \cdot \frac{d\omega}{ds}$$

$$= R'\left(\frac{(R\kappa - 1)^2 + R'^2 + (R\kappa - 1)^2 - R'^2}{(R\kappa - 1)^2 + R'^2}\right) + R\cos\omega \cdot \frac{d\omega}{ds}$$

$$= (R\kappa - 1)\left(\frac{2R'(R\kappa - 1)}{(R\kappa - 1)^2 + R'^2}\right) + R\cos\omega \cdot \frac{d\omega}{ds}$$

$$= (R\kappa - 1)\cos\omega + R\cos\omega \cdot \frac{d\omega}{ds}$$

$$= \cos\omega\left[R\kappa - 1 + R\frac{d\omega}{ds}\right] = \cos\omega \cdot 0 = 0 \qquad \blacksquare$$

------

## 7. Observations and More Examples

If the radiant S is at infinity, then the radius of the focal circle $C_s$ is R(s) = $\frac{1}{4\kappa(s)}$, where $\kappa(s)$ is the curvature of the reflective curve $\boldsymbol{\alpha}(s)$. In this case, the envelope $\boldsymbol{\beta}(s)$ is given by Corollary 3. Consequently, the same family of focal circles $C_s$ and their envelope $\boldsymbol{\beta}(s)$ can be used to generate caustics for all radiants at infinity. As the focal circles roll on $\boldsymbol{\beta}(s)$, different points on $C_s$ will simultaneously trace out the caustics resulting from different radiants at infinity. In particular, for each s, the point of contact P between $C_s$ and $\boldsymbol{\beta}(s)$ is generating one of these caustics. Theorem 1 tells us that this caustic has a cusp at P. This gives us an alternate description of $\boldsymbol{\beta}(s)$ as the locus of cusps of caustics resulting from radiants at infinity. As the direction of the radiant at infinity is varied, the cusps of caustics will slide along $\boldsymbol{\beta}(s)$. This suggests a sundial can be constructed utilizing this property. Reflected light at $\boldsymbol{\alpha}(s)$ will form a cusp if the incident light ray approaches along the axis of aberrancy.

As a slight disclaimer, it should pointed out that in all of the following examples, the caustics that are described explicitly by an equation or otherwise have been described elsewhere using other methods. Our purpose is to use the Theorems and Corollaries of this paper to describe these caustics.

**Example 3** (continued): We showed that for the parabola with radiants at infinity, the envelope $\boldsymbol{\beta}(s)$ collapses to a single point, the focus of the parabola. The rolling of the focal circles along $\boldsymbol{\beta}(s)$ simply amounts to pivoting around the focus. See [1] for the animation of this example.

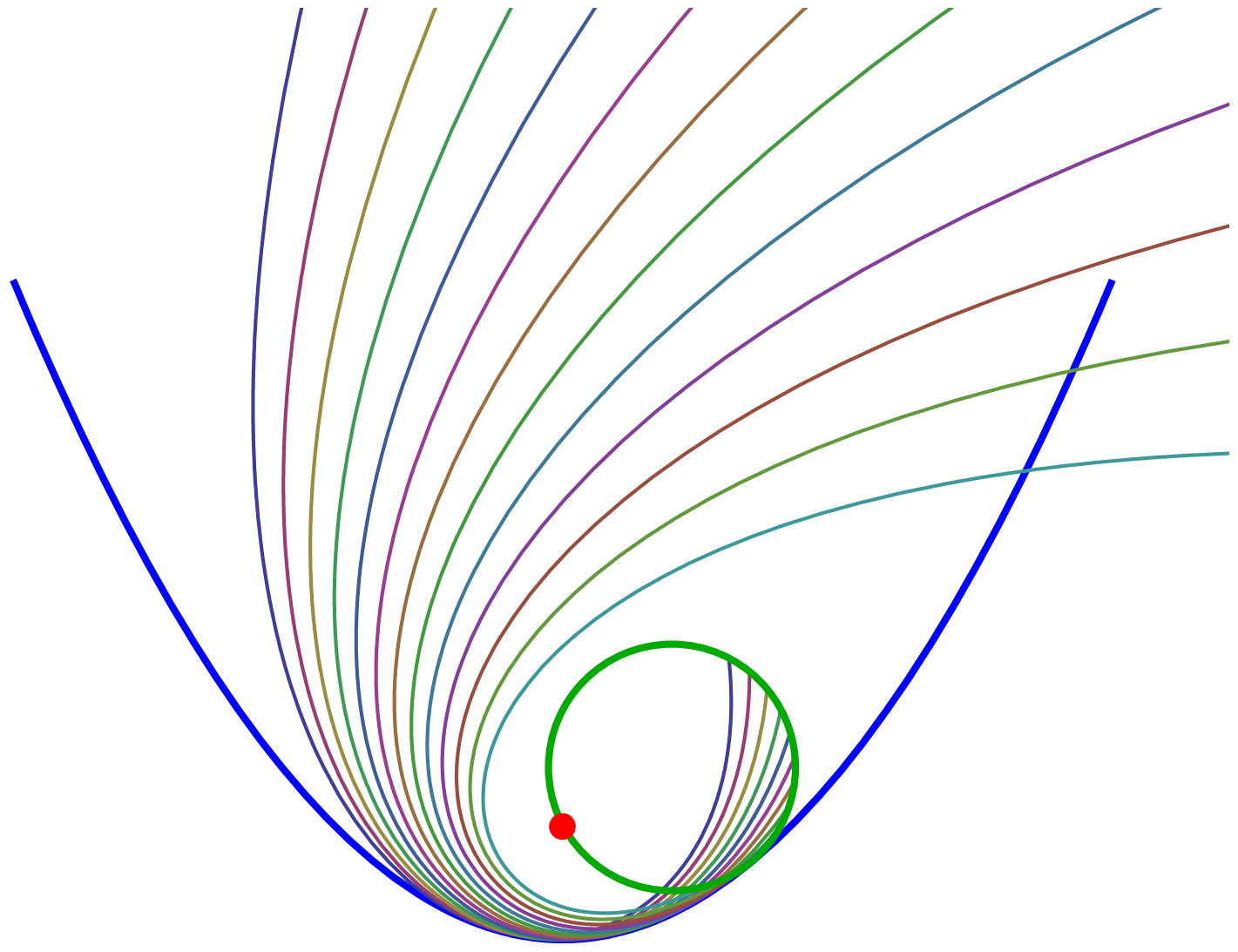

**Figure 13** Generating multiple caustics from radiants at infinity

Figure 13 illustrates the simultaneous creation of caustics from thirteen radiants at infinity coming in from the directions -.3 to .9 radians in .1 radian increments. The focal circles slide tangent to the parabola as they pivot around the focus. The thirteen points generating the caustics maintain fixed angular separations on the focal circles. If the radiant is in the direction $\frac{\pi}{2}$ parallel to the axis of the parabola, then the generating point of the caustic is the pivot point, the focus, thus confirming the entire caustic collapses to the focus. From this we observe the well known property of the parabola that for radiants at infinity, there are no cusps on the caustics (unless the direction is parallel to the axis). It seems that parabolas are the unique curves with this property. However, if the radiant is not at infinity, a cusp can be formed.

For the parabola $\boldsymbol{\alpha}(t) =< t, t^2 >$, we will find a parameterization for the caustic generated by the point $P_0$ on the focal circle $C_t$ opposite the focus F$(0, \frac{1}{4})$ resulting from reflected light coming in parallel to the x-axis before reflection. The focal circle $C_t$ has radius $\frac{1}{4\kappa} = \frac{1}{8}(4t^2+1)^{3/2}$ and center $c(t) =< t, t^2 > + \frac{1}{8}(4t^2+1)^{3/2} \frac{<-2t,1>}{(4t^2+1)^{1/2}}$

$= \frac{1}{8} < 6t - 8t^3, 12t^2 + 1 >$. The caustic is then $\boldsymbol{E}(t) =< 0, \frac{1}{4} > +2\left(c(t) - < 0, \frac{1}{4} >\right) = \frac{1}{2} < 3t - 4t^3, 6t^2 >$.

Eliminating the parameter yields $108x^2 = y(4y-9)^2$. This curve is a version of a Tschirnhausen Cubic curve [18].

Now consider the caustic generated by the point $P_\delta$ on $C_t$ pictured in Figure 14. $P_0F$ is a diameter of $C_t$ and the chord $P_\delta F$ makes the angle $\delta$ with $P_0F$. The point $P_0$ generates the caustic considered above. The dynamics make it clear that the caustic generated by $P_\delta$ is obtained from the caustic generated by $P_0$ by rotating about F through the angle $\delta$ and contracting toward F by a factor $\cos\delta$. Hence all such caustics are essentially the same curve.

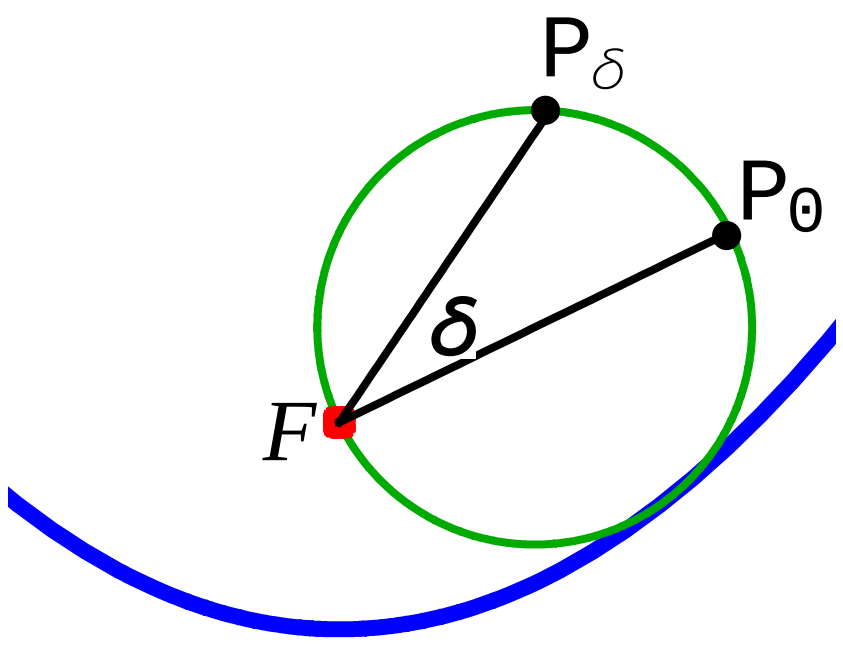

**Figure 14** Points generating two caustics

**Example 4:** Consider the deltoid $\boldsymbol{\alpha}(t) = < 2\cos(t) + \cos(2t), 2\sin(t) - \sin(2t) >$ with radiants at infinity. Using Corollary 3, we calculate $\boldsymbol{\beta}(t) =< 4\cos(t) - \cos(4t), 4\sin(t) - \sin(4t) >$. Hence, $\boldsymbol{\alpha}(t)$ and $\boldsymbol{\beta}(t)$ are a hypocycloid and epicycloid traced out by points on unit circles as they roll around the inside and outside of a circle of radius 3, respectively. The focal circles have their centers on this circle of radius 3. Remarkably, the caustic itself is an astroid that is circumscribed about the deltoid and inscribed in the hypocycloid $\boldsymbol{\beta}$ for every radiant at infinity. See [5] and Figure 15.

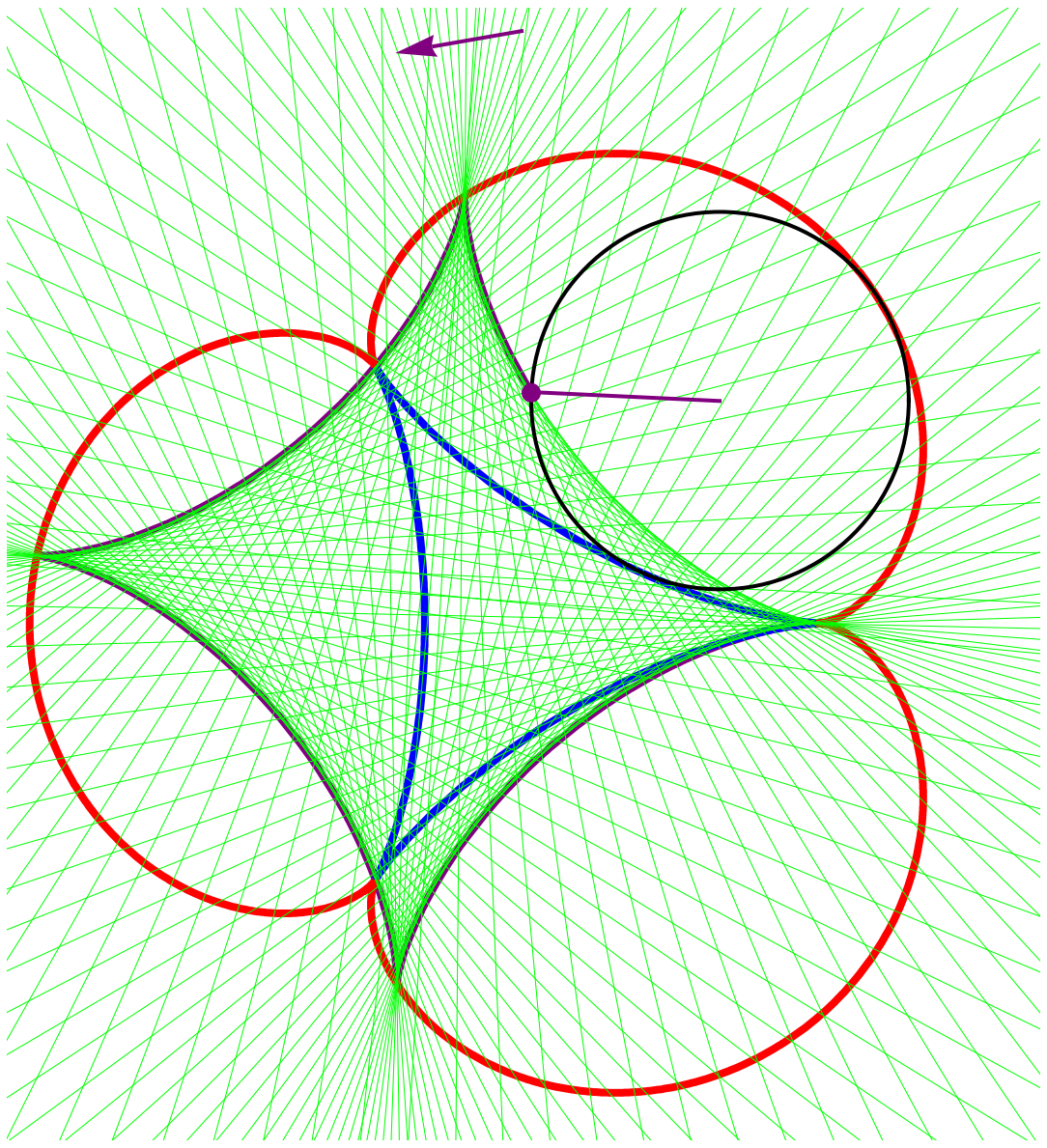

**Figure 15** Tracing the astroidal caustic of the deltoid

See [1] for a stunning animation of this example. By varying the angle of the radiant at infinity, the astroid rotates with its 4 cusps sliding around $\boldsymbol{\beta}$ while at the same time maintaining contact with the 3 cusps of the deltoid. This effect is also animated at [1].

The next two examples use a circular mirror. For detailed calculations for all caustics from a circular mirror, see Cayley [11] where he provides equations for the caustics along with careful descriptions of their significant features (also see [4]).

**Example 5:** First consider the circular mirror with internal reflection from a radiant S that is on the circle itself. For the mirror we will use the unit circle U, $\boldsymbol{\alpha}(s)$ = <cos $s$, sin $s$>. At each point of the circle we have two more circles tangent to U. One is the circle that passes through the radiant S, and the other is the focal circle $C_s$ that passes through the point F on the caustic resulting from reflection at $\boldsymbol{\alpha}(s)$. The radii of these two circles are related by Theorem 2. Since S is itself on U, it is apparent that the circle passing through S is U itself, independent of the particular point of reflection $\boldsymbol{\alpha}(s)$. This circle has diameter $d_1$ = 2 and constant curvature $\kappa$ = 1. Theorem 2 then says each focal circle $C_s$ has diameter $d_2 = \frac{2}{3}$, and from this we deduce that the envelope $\boldsymbol{\beta}(s)$ is a circle of diameter $\frac{2}{3}$ centered at the origin. Hence, the caustic is another epicycloid traced out by a point on the circle of radius $\frac{1}{3}$ rolling around the outside of a circle of radius $\frac{1}{3}$. See Figures 16 and 17 and Photo 1. To see the animation of this example see [2].

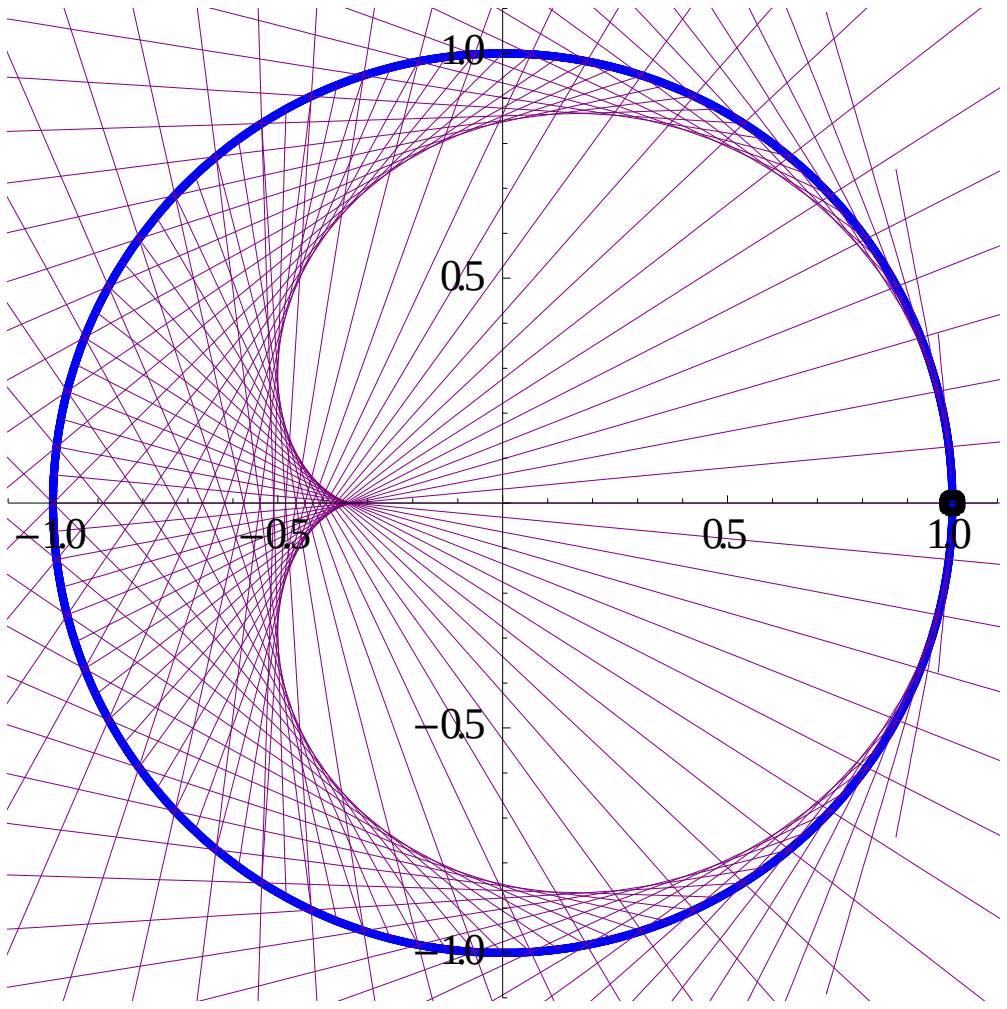

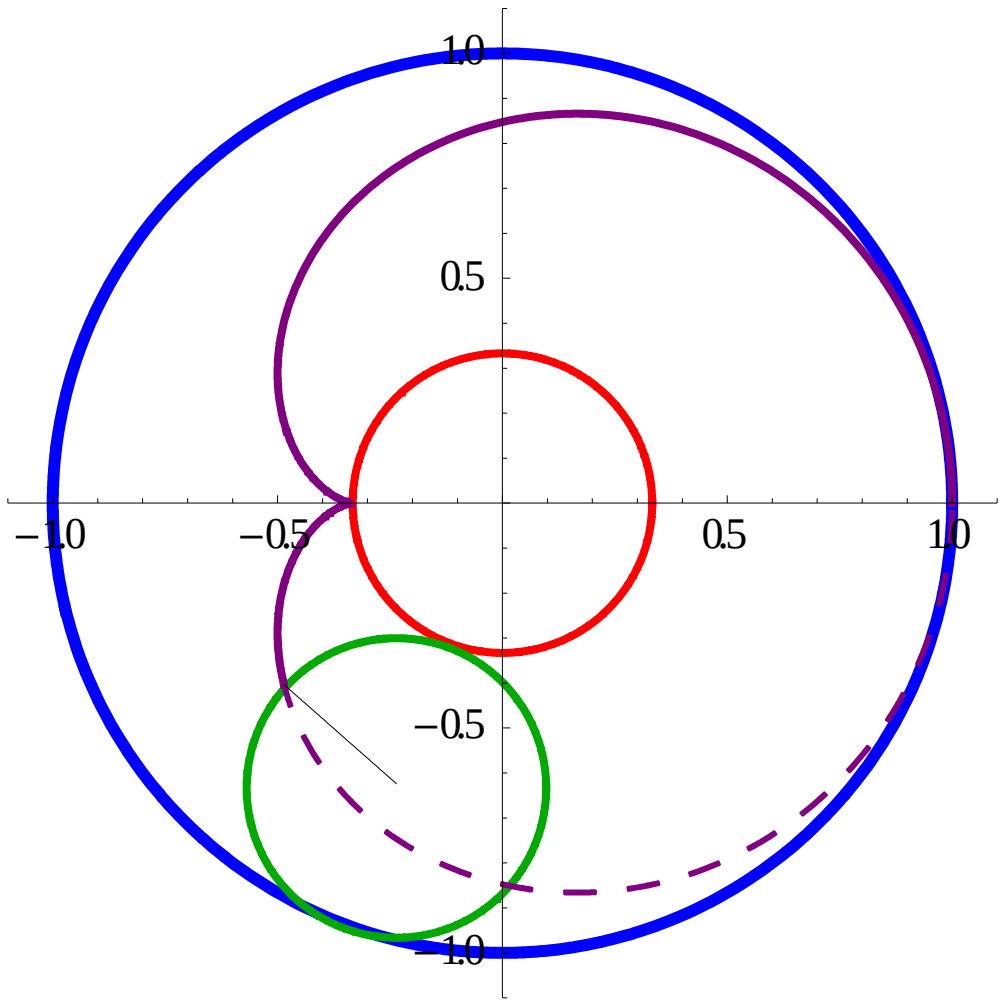

**Figure 16** Reflection from radiant on circular mirror

**Figure 17** Tracing the epicycloidal caustic

**Example 6:** Next, consider the circular mirror with a radiant on the interior. Utilizing the unit circle U again, we place the radiant at (c,0), where $0 < c < \frac{1}{2}$. The resulting caustic is a closed curve with four cusps. When $\frac{1}{2} < c < 1$, the caustic has two components . Figure 18 shows the caustic when c = .25. Also, see Photo 2. The envelope $\boldsymbol{\beta}(s)$ also has four cusps.

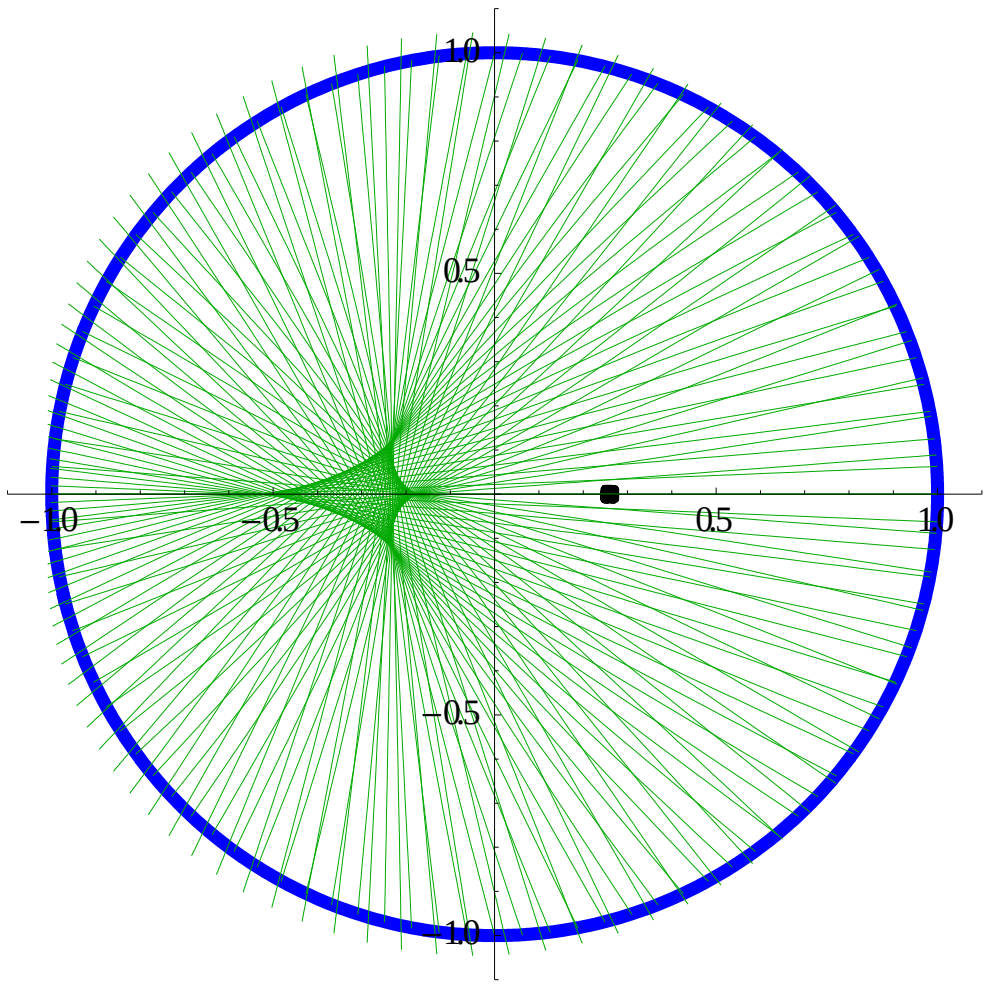

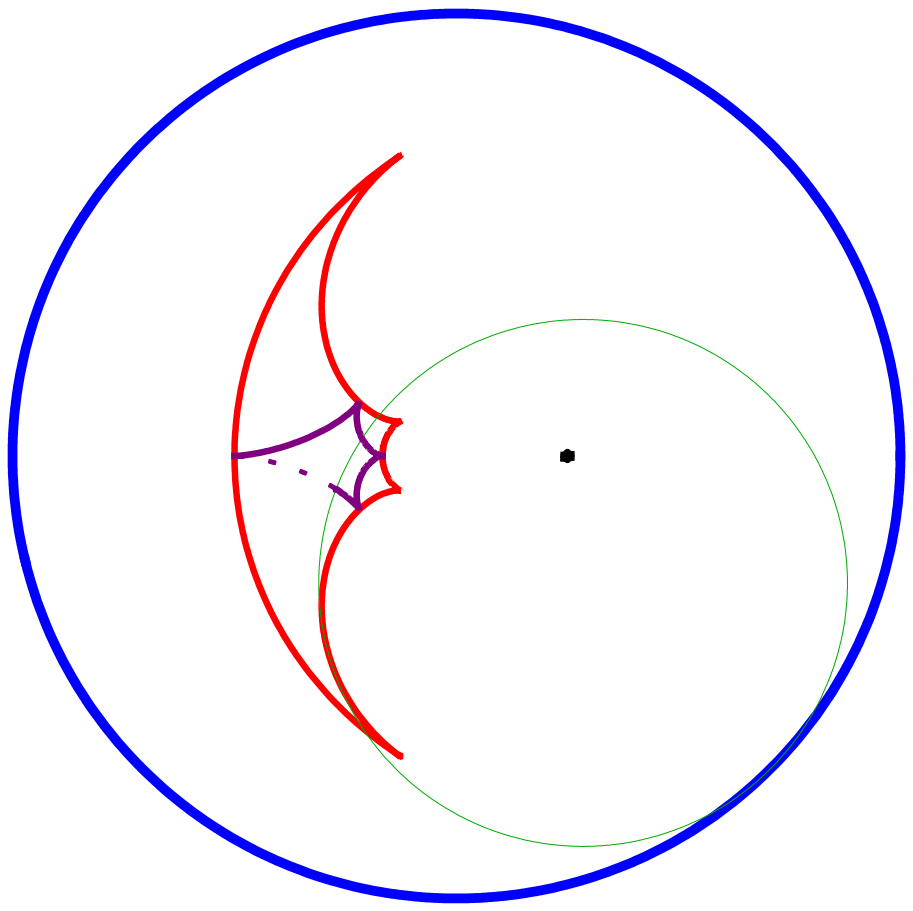

**Figure 18** Circular mirror with interior radiant

**Figure 19** Tracing the caustic

Figure 19 shows the focal circle rolling on $\boldsymbol{\beta}(s)$ to generate the caustic. The focal circles roll seamlessly and continuously from one smooth section of $\boldsymbol{\beta}(s)$ across the cusps to the next smooth section without creating any noticeable distinguishing feature of the caustic.

## 8. Final Remarks

We finish with some miscellaneous observations about features of caustics. In the following, when we speak of the radiant moving, we really mean the apparent location of the radiant shifts as we move along a section of the reflective curve.

When the radiant is sufficiently close to the reflective curve, the caustic can have two or more components. Recalling the observations following Theorem 2, at a particular point of the reflective curve, if the radiant is on the concave side and outside the discriminant circle, then the corresponding point of the caustic will also be on the concave side outside the discriminant circle. As the radiant moves closer to and then intersects the discriminant circle, the point of the caustic migrates out to infinity. As the radiant then passes to the inside of the discriminant circle, the caustic returns from infinity, but in the opposite direction on the convex side of the curve. The reflected line at the moment the radiant encounters a discriminant circle will serve as an asymptote for the caustic as it goes to and returns from infinity. For example, Figure 5 illustrates an elliptical $\boldsymbol{\alpha}$ and the corresponding $\boldsymbol{\beta}$ for radiants at infinity (i.e., the second envelope resulting from the family of discriminant circles). For any radiant either exterior to the ellipse or interior to $\boldsymbol{\beta}$ the caustic will be a single component curve. If the radiant is between the ellipse and $\boldsymbol{\beta}$, then the caustic will have two components. In Figure 1, the purple caustic has a second component which in not seen since it is beyond the window of the graphic.

Now suppose along a section of the reflective curve the radiant passes from the concave side to the convex side. Let P be the point on the curve whose tangent line contains the radiant. As the radiant approaches the tangent line, the circle tangent to the reflective curve at P that contains the radiant will have diameter that goes to infinity. Consequently, the focal circle that contains the corresponding point of the caustic will approach the discriminant circle at P. But since the angle of incidence is the angle of reflection, the point of the caustic must approach P. As the radiant passes to the other side of the tangent line (the convex side at P), the caustic emerges from P on the inside of the discriminant circles.

**Acknowledgements.** I would like to thank Todd Will for creating the demonstrations at [1, 2] that superbly illustrate the Theorems of this paper. Also, I am grateful to the anonymous referees and the editor for the careful reading of the article and their many helpful improvements. Also, I thank Robert Allen who graciously volunteered to help with the final preparation of the paper.

**Jeffrey A. Boyle** received his Ph.D. from the University of Iowa in 1984, taught at Michigan State University, and from 1988 to present, at the University of Wisconsin-La Crosse. He is an avid cross country skier, cyclist, runner, and most importantly, dog walker.

*Department of Mathematics, University of Wisconsin-La Crosse, La Crosse, WI 54601*
*jboyle@uwlax.edu*